\documentclass{article}%
\usepackage{amsmath}
\usepackage{color}
\usepackage{amsfonts}
\usepackage{amssymb}
\usepackage{graphicx}%
\usepackage{caption}
\usepackage{subcaption}
\usepackage{float}
\providecommand{\U}[1]{\protect\rule{.1in}{.1in}}
\newtheorem{theorem}{Theorem}

\newtheorem{corollary}[theorem]{Corollary}

\newtheorem{lemma}[theorem]{Lemma}

\newtheorem{remark}[theorem]{Remark}

\newenvironment{proof}[1][Proof]{\noindent\textbf{#1.} }{\ \rule{0.5em}{0.5em}}
\begin{document}

\title{On a retarded stochastic system with discrete diffusion modeling life tables}
\author{Tom\'{a}s Caraballo$^{1}$, Francisco Morillas$^{2}$, Jos\'{e} Valero$^{3}$\\$^{1}$ {\small Universidad de Sevilla, Departamento de Ecuaciones
Diferenciales y An\'{a}lisis Num\'{e}rico}\\{\small Apdo. de Correos 1160, 41080-Sevilla Spain}\\{\small E.mail:\ caraball@us.es}\\$^{2}${\small Universitat de Val\`{e}ncia, Departament d'Economia Aplicada,
Facultat d'Economia, }\\{\small Campus dels Tarongers s/n, 46022-Val\`{e}ncia, Spain. }\\{\small E.mail: Francisco.Morillas@uv.es}\\$^{3}${\small Universidad Miguel Hern\'{a}ndez de Elche, Centro de
Investigaci\'{o}n Operativa} \\{\small Avda. Universidad s/n, Elche (Alicante), 03202, Spain. }\\{\small E.mail: jvalero@umh.es}}
\date{}
\maketitle

\begin{abstract}
This work proposes a method for modeling and forecasting mortality rates. It
constitutes an improvement over previous studies by incorporating both the
historical evolution of the mortality phenomenon and its random behavior. In
the first part, we introduce the model and analyze mathematical properties
such as the existence of solutions and their asymptotic behavior. In the
second part, we apply this model to forecast mortality rates in Spain, showing
that it yields better results than classical methods.

\end{abstract}

\bigskip

\section{Introduction}

In actuarial or demographic sciences, life tables are very useful to study
biometric functions such as the probability of survival or death, therefore
they are crucial to calculate the insurance premium. In our previous papers
\cite{MoVa14}, \cite{MoVa01}, \cite{CarMoVa02}, we have shown that a system of
ordinary differential equations with nonlocal discrete diffusion is
appropriate to model dynamical life tables. In \cite{MoVa14}, we implemented
the deterministic model
\begin{align}
\frac{d}{dt}u_{i}\left(  t\right)   &  =\sum_{r\in D}j_{i-r}u_{r}\left(
t\right)  -u_{i}(t)+\sum_{r\in\mathbb{Z}\backslash D}j_{i-r}g_{r}\left(
t\right)  \text{, }i\in D\text{, }t>0,\label{NoLocalEq}\\
u_{i}\left(  0\right)   &  =u_{i}^{0}\text{, }i\in D\text{,}\nonumber
\end{align}
where $D=\left\{  m_{1},\ldots,m_{2}\right\}  $, $m_{1}<m_{2}$, $m_{i}%
\in\mathbb{Z}$, and $j:\mathbb{Z}\rightarrow\mathbb{R}^{+}$. Using the
observed data published by Spain's National Institute of Statistics, we
compared the results of our model with those obtained through classical
techniques. Our numerical simulations indicate that, in the short
run---specifically for predictions within three years---we achieved
improvements in certain indicators of goodness and smoothness. Since some
memory effects are present on life tables, we considered in \cite{MoVa01} a
modification of model (\ref{NoLocalEq}) in which we added some delay terms.
Namely, we studied the problem:%
\begin{align}
\dfrac{d}{dt}u_{i}\left(  t\right)   &  =\sum_{r\in D}\int_{-h}^{0}%
j_{i-r}u^{r}(t+s)\alpha_{i}(s)d\mu(s)-u_{i}(t)\label{Delay}\\
&  +\sum_{r\in\mathbb{Z}\backslash D}\int_{-h}^{0}j_{i-r}g_{r}\left(
t+s\right)  \alpha_{i}(s)d\mu(s)\text{, }i\in D,\ t>0,\nonumber\\
u^{i}\left(  \tau+s\right)   &  \equiv\phi^{i}\left(  s\right)  \text{, }i\in
D\text{, }s\in\lbrack-h,0],\nonumber
\end{align}
where $\alpha_{i}:[-s,0]\rightarrow\mathbb{R}^{+}$ and $d\mu(s)=\xi\left(
s\right)  ds$ being $\xi\left(  \text{\textperiodcentered}\right)  $ a
probability density. In \cite{MoVa01} it was shown that with this model the
prediction horizon can be extended up to $8$ years. In addition, it gives
coherent values, in magnitude, when comparing it with other classical
techniques such as the Lee-Carter model, up to $18$ years.

Despite the good results provided by these deterministic models, in the real
world there is always a\ certain level of noise, which either can be intrinsic
to the model or can appear due to the presence of errors in the observed data.
This is why in \cite{CarMoVa02} we considered the stochastic version of model
(\ref{NoLocalEq}) given by%
\begin{align}
\frac{d}{dt}u_{i}\left(  t\right)   &  =\sum_{r\in D}j_{i-r}u_{r}\left(
t\right)  -u_{i}(t)+\sum_{r\in\mathbb{Z}\backslash D}j_{i-r}g_{r}\left(
t\right)  +b\sigma_{i}(u_{i}(t))\frac{dw_{i}}{dt}\text{, }i\in D\text{,
}t>0,\label{NoLocalStoch}\\
u_{i}\left(  0\right)   &  =u_{i}^{0}\text{, }i\in D\text{,}\nonumber
\end{align}
where $w_{i}\left(  t\right)  $ are independent Brownian motions and $b>0$ is
the intensity of the white noise, and two specific type of noises were
considered: 1) $\sigma_{i}(v)=v$ (linear case); 2) $\sigma_{i}(v)=v(1-v)$. The
choice of the noise in the second case is motivated by the fact that we are
interested in studying variables like the probability of death, which take
values in the interval $[0,1]$. Although the results given by the numerical
simulations are fine from the qualitative point of view, it was necessary to
make a correction of the estimates by using the average annual improvement
rate. However, this correction is not equally adequate for all ages. Thus, in
order to take into account in the model appropriate correction rates for each
age, we need to introduce delay terms in the equation as in (\ref{Delay}).
Therefore, we will study in this paper the following stochastic system of
differential equations with delay:%
\begin{equation}
\left\{
\begin{array}
[c]{c}%
\dfrac{d}{dt}u^{i}\left(  t\right)  =\left[  J(t,u_{t})\right]  _{i}%
-u^{i}(t)+b\sigma_{i}(u_{i}(t))\dfrac{dw_{i}}{dt}\text{, }i\in D\text{,
}t>\tau,\\
u^{i}\left(  \tau+s\right)  \equiv\phi^{i}\left(  s\right)  \text{, }i\in
D\text{, }s\in\lbrack-h,0],
\end{array}
\right.  \label{DelayStoch}%
\end{equation}
where $D=\{m_{1},m_{1}+1,...,m_{2}\}$, $-\infty<m_{1}<m_{2}<\infty$,
$m=m_{2}-m_{1}+1,$ $\tau$ is the initial moment of time, $h>0$, $w_{i}\left(
t\right)  $ are independent Brownian motions, $b\geq0$ is the intensity of the
white noise, $\sigma_{i}:\mathbb{R}\rightarrow\mathbb{R},$ and $J:\mathbb{R}%
^{+}\times C([-h,0],\mathbb{R}^{m})\rightarrow\mathbb{R}^{m}$ is the
non-autonomous convolution operator defined by
\[
\left[  J(t,u_{t})\right]  _{i}=\sum_{r\in D}\int_{-h}^{0}j_{i-r}%
u^{r}(t+s)\alpha_{i}(s)d\mu(s)+\sum_{r\in\mathbb{Z}\backslash D}\int_{-h}%
^{0}j_{i-r}g_{r}\left(  t+s\right)  \alpha_{i}(s)d\mu(s)\text{, if }i\in
D\text{,}%
\]
where $u\left(  t\right)  =\left(  u^{i}\left(  s\right)  \right)  _{i\in D},$
$u_{t}$ is the segment of solution defined by $u_{t}(s)=u\left(  t+s\right)  $
for $s\in\lbrack-h,0]$, $j:\mathbb{Z}\rightarrow\mathbb{R}$, $g:\mathbb{R}%
\times\mathbb{Z}\backslash D\rightarrow\mathbb{R}$ and $d\mu(s)=\xi\left(
s\right)  ds$ being $\xi\left(  \text{\textperiodcentered}\right)  $ a
probability density. We define the Banach space%
\[
l_{2}^{\infty}=\{\left(  u_{i}\right)  _{i\in\mathbb{Z}\backslash D}%
:\sup_{i\in\mathbb{Z}\backslash D}\left\vert u_{i}\right\vert <\infty\}
\]
with the norm $\left\Vert u\right\Vert _{l_{2}^{\infty}}=\sup_{i\in
\mathbb{Z}\backslash D}\left\vert u_{i}\right\vert $ and assume the following conditions:

\begin{itemize}
\item[$\left(  H1\right)  $] $j_{k}\geq0$ for all $k\in\mathbb{Z}$.

\item[$\left(  H2\right)  $] $\sum_{i\in\mathbb{Z}}j_{i}=1$.

\item[$\left(  H3\right)  $] $g\in C(\mathbb{R},l_{2}^{\infty}).$

\item[$\left(  H4\right)  $] $\alpha_{i}\in C([-h,0],\mathbb{R}),\ \alpha
_{i}\left(  s\right)  \geq0$ for $i\in D,\ s\in\lbrack-h,0]$.
\end{itemize}

As for problem (\ref{NoLocalStoch}) we will consider the noises $\sigma
_{i}(v)=v$ (linear case) and $\sigma_{i}(v)=v(1-v).$

In the first part of the paper, we study several properties of the solutions
to problem (\ref{DelayStoch}). In the case where the noise is linear, we
establish the existence of a unique globally defined positive solution if the
initial condition is positive. When the noise is non-linear with $\sigma
_{i}(v)=v(1-v)$, we prove that if the initial condition belongs to $\left(
0,1\right)  $, then the solution remains in this interval for any future
moment of time. Finally, we analyze the asymptotic behaviour of the solutions
as time goes to $+\infty$, showing under certain assumptions that, for large
enough time, the solution belongs to a neighbourhood of the unique fixed point
of the deterministic system.

In the second part of the paper, we apply model (\ref{DelayStoch}) to forecast
age-specific mortality rates in Spain. Using observed data from 2008 to 2018,
we perform numerical simulations to generate multiple realizations of
predicted mortality values for the period 2019--2023. Based on these
realizations, we construct confidence intervals and calculate several error
indicators, comparing the results with those obtained from classical
techniques such as the Lee-Carter and Renshaw-Haberman models. Our model
achieved the best results for all years within the validation period
(2019-2023). Thus, we conclude that our method should be regarded as a
promising alternative to classical models.

\section{Properties of solutions}

In this section, we shall obtain some properties of the solutions to problem
(\ref{DelayStoch}).

\section{The linear case}

We will first consider a standard linear noise, that is, we study the system%
\begin{align}
\frac{d}{dt}u_{i}\left(  t\right)   &  =\left[  J(t,u_{t})\right]  _{i}%
-u_{i}(t)+bu_{i}(t)\frac{dw_{i}}{dt}\text{, }i\in D\text{, }t>\tau
,\label{Linear}\\
u_{i}\left(  \tau+s\right)   &  \equiv\phi^{i}\left(  s\right)  \text{, }i\in
D\text{, }s\in\lbrack-h,0]\text{.}\nonumber
\end{align}

Denote $\mathbb{R}_{+}^{m}=\{v\in\mathbb{R}^{m}:v_{j}>0$ for all $j\}$ and
$\overline{\alpha}_{i}=\int_{-h}^{0}\alpha_{i}(s)d\mu(s)$. Our aim is to
establish the existence of global positive solutions.

\begin{lemma}
\label{Positive}Assume that $g_{r}(t)\geq0,$ for all $r$ and $t,$ and that
$\overline{\alpha}_{i}\leq1$, for all $i$. Then for any $\phi\in
C([-h,0],\mathbb{R}_{+}^{m})$ there exists a unique globally defined solution
$u\left(  \text{\textperiodcentered}\right)  $ such that $u\left(  t\right)
\in\mathbb{R}_{+}^{m}$ almost sure for $t\geq\tau.$
\end{lemma}

\begin{proof}
The existence of a unique local solution to problem (\ref{Linear}) follows
from standard results for functional stochastic differential equations
governed by locally Lipschitz functions \cite[Theorem 2.8, P. 154]{MaoBook}.
Given that any solution $u\left(  \text{\textperiodcentered}\right)  $ is
defined in the maximal interval $[0,\tau_{e})$, we need to prove that
$\tau_{e}=+\infty$ and that $u\left(  t\right)  \in\mathbb{R}_{+}^{m}$ for all
$t\geq0$ a.s.

We choose $k_{0}>0$ such that%
\[
\frac{1}{k_{0}}<\min_{s\in\lbrack-h,0]}\left\vert \phi_{i}(s)\right\vert
\leq\max_{s\in\lbrack-h,0]}\left\vert \phi_{i}(s)\right\vert <k_{0}\text{ for
all }i\in D.
\]
For each $k\geq k_{0}$ we define the stopping time%
\[
\tau_{k}=\inf\{t\in\lbrack\tau,\tau_{e}):u_{i}\left(  t\right)  \not \in
(\frac{1}{k},k)\text{ for some }i\in D\}.
\]
This sequence is increasing as $k\nearrow+\infty$. If $\tau_{\infty}%
=\lim_{k\rightarrow+\infty}\tau_{k}=+\infty$ a.s., then $\tau_{e}=+\infty$ and
$u\left(  t\right)  \in\mathbb{R}_{+}^{m}$ for $t\geq0$ almost surely, proving
the assertion.

If $\lim_{k\rightarrow+\infty}\tau_{k}\not =+\infty$ a.s, there would exist
$T,\varepsilon>0$ such that
\[
P(\tau_{\infty}\leq T)>\varepsilon,
\]
and then there would be $k_{1}\geq k_{0}$ for which%
\[
P(\tau_{k}\leq T)\geq\varepsilon\text{ for }k\geq k_{1}.
\]
Further, we consider the $C^{2}$-function $V:\mathbb{R}_{+}^{m}\rightarrow
\mathbb{R}_{+}^{1}$ given by%
\[
V\left(  u\right)  =\sum_{i\in D}(u_{i}-1-\log(u_{i})).
\]

Let $\tau\leq t\leq\tau_{k}\wedge T.$ Then $u\left(  t\right)  \in
\mathbb{R}_{+}^{m}$ and by It\^{o}'s formula we have%
\begin{align}
&  dV(u(t))\label{EqV}\\
&  =\sum_{i\in D}\left(  1-\frac{1}{u_{i}\left(  t\right)  }\right)  \left(
\sum_{r\in D}\int_{-h}^{0}j_{i-r}u_{r}(t+s)\alpha_{i}(s)d\mu(s)+\sum
_{r\in\mathbb{Z}\backslash D}\int_{-h}^{0}j_{i-r}g_{r}\left(  t+s\right)
\alpha_{i}(s)d\mu(s)-u_{i}(t)\right)  dt\nonumber\\
&  +\sum_{i\in D}\frac{1}{2}b^{2}dt+\sum_{i\in D}\left(  1-\frac{1}%
{u_{i}\left(  t\right)  }\right)  bu_{i}(t)dw_{i}(t)\nonumber\\
&  =I(t)dt+\sum_{i\in D}I_{i}(t)dw_{i}(t),\nonumber
\end{align}
where $I_{i}(t)=\left(  1-\frac{1}{u_{i}\left(  t\right)  }\right)  bu_{i}%
(t)$. By using $\left(  H1\right)  -\left(  H4\right)  $, $u\left(  t\right)
\in\mathbb{R}_{+}^{m}$ and $g_{r}\left(  t\right)  \geq0$ the first term is
estimated by:%
\begin{align*}
I(t)  &  \leq\sum_{i\in D}\sum_{r\in D}\int_{-h}^{0}j_{i-r}u_{r}%
(t+s)\alpha_{i}(s)d\mu(s)-\sum_{i\in D}u_{i}(t)+m+\frac{m}{2}b^{2}\\
&  +\sum_{i\in D}\sum_{r\in\mathbb{Z}\backslash D}\int_{-h}^{0}j_{i-r}%
g_{r}\left(  t+s\right)  \alpha_{i}(s)d\mu(s)\\
&  -\sum_{i\in D}\frac{1}{u_{i}\left(  t\right)  }\left(  \sum_{r\in D}%
\int_{-h}^{0}j_{i-r}u_{r}(t+s)\alpha_{i}(s)d\mu(s)+\sum_{r\in\mathbb{Z}%
\backslash D}\int_{-h}^{0}j_{i-r}g_{r}\left(  t+s\right)  \alpha_{i}%
(s)d\mu(s)\right) \\
&  \leq K_{T}+\sum_{i\in D}\sum_{r\in D}\int_{-h}^{0}j_{i-r}u_{r}%
(t+s)\alpha_{i}(s)d\mu(s)-\sum_{i\in D}u_{i}(t),
\end{align*}
where we have used that%
\[
\sum_{i\in D}\frac{1}{u_{i}\left(  t\right)  }\left(  \sum_{r\in D}\int%
_{-h}^{0}j_{i-r}u_{r}(t+s)\alpha_{i}(s)d\mu(s)+\sum_{r\in\mathbb{Z}\backslash
D}\int_{-h}^{0}j_{i-r}g_{r}\left(  t+s\right)  \alpha_{i}(s)d\mu(s)\right)
\geq0,
\]%
\[
\sum_{i\in D}\sum_{r\in\mathbb{Z}\backslash D}\int_{-h}^{0}j_{i-r}g_{r}\left(
t+s\right)  \alpha_{i}(s)d\mu(s)\leq mC_{T}.
\]

Integrating in (\ref{EqV}) over $\left(  \tau,\tau_{k}\wedge T\right)  $ and
taking expectations we obtain that%
\begin{align*}
0  &  \leq\mathbb{E}V(u(\tau_{k}\wedge T))\\
&  \leq V(\phi(0))+\mathbb{E}\int_{\tau}^{\tau_{k}\wedge T}K_{T}%
dt+\mathbb{E}\int_{\tau}^{\tau_{k}\wedge T}\sum_{i\in D}\left(  1-\frac
{1}{u_{i}\left(  t\right)  }\right)  bu_{i}(t)dw_{i}(t)\\
&  +\mathbb{E}\sum_{i\in D}\sum_{r\in D}\int_{-h}^{0}\int_{\tau}^{\tau
_{k}\wedge T}j_{i-r}u_{r}(t+s)\alpha_{i}(s)dtd\mu(s)-\mathbb{E}\sum_{i\in
D}\int_{\tau}^{\tau_{k}\wedge T}u_{i}(t)dt.
\end{align*}
Using $\overline{\alpha}_{i}\leq1$ we have%
\begin{align*}
&  \sum_{i\in D}\sum_{r\in D}\int_{-h}^{0}\int_{\tau}^{\tau_{k}\wedge
T}j_{i-r}u_{r}(t+s)\alpha_{i}(s)dtd\mu(s)\\
&  =\sum_{i\in D}\sum_{r\in D}\int_{-h}^{0}j_{i-r}\alpha_{i}(s)\int_{\tau
}^{\tau_{k}\wedge T}u_{r}(t+s)dtd\mu(s)\\
&  \leq\sum_{r\in D}\sum_{i\in D}j_{i-r}\int_{-h}^{0}\alpha_{i}(s)d\mu
(s)\int_{\tau-h}^{\tau_{k}\wedge T}u_{r}(t)dt\\
&  \leq\sum_{r\in D}\int_{\tau-h}^{\tau_{k}\wedge T}u_{r}(t)dt=\sum_{r\in
D}\int_{\tau}^{\tau_{k}\wedge T}u_{r}(t)dt+\sum_{r\in D}\int_{\tau-h}^{\tau
}u_{r}(t)dt.
\end{align*}
Hence,%
\begin{align*}
0  &  \leq\mathbb{E}V(u(\tau_{k}\wedge T))\\
&  \leq V(\phi(0))+K_{T}\mathbb{E(}\tau_{k}\wedge T)+\sum_{r\in D}\int%
_{\tau-h}^{\tau}\phi(t)dt\\
&  \leq V(\phi(0))+K_{T}T+K_{\phi}.
\end{align*}

Let $\Omega_{k}=\{\omega:\tau_{k}\leq T\}$, which satisfies $P(\Omega_{k}%
)\geq\varepsilon$ for $k\geq k_{1}$. For any $\omega\in\Omega_{k}$ there is
$i\in D$ such that either $u_{i}\left(  \tau_{k},\omega\right)  =k$ or
$u_{i}\left(  \tau_{k},\omega\right)  =1/k$, which implies that
\[
V(u(\tau_{k}\wedge T,\omega))\geq(k-1-\log(k))\wedge(\frac{1}{k}-1+\log(k)).
\]
Hence,%
\[
V(\phi(0))+K_{T}T+K_{\phi}\geq\mathbb{E}(1_{\Omega_{k}}V(u(\tau_{k}\wedge
T)))\geq\varepsilon((k-1-\log(k))\wedge(\frac{1}{k}-1+\log(k)))=\varepsilon
R(k),
\]
where $1_{A}$ stands for the indicator function of the set $A$. Passing to the
limit as $k\rightarrow+\infty$ we obtain a contradiction as $R(k)\rightarrow
+\infty$.
\end{proof}

\bigskip As a consequence, the following result follows.

\begin{corollary}
Let $\phi^{1},\phi^{2}\in C([-h,0],\mathbb{R}_{+}^{m})$ be two initial
conditions satisfying $\phi_{i}^{1}(s)>\phi_{i}^{2}(s)$ for any $i\in D$,
$s\in\lbrack-h,0]$. Also, $g^{1},g^{1}\in C([0,+\infty),l_{2}^{\infty})$ are
such that $g_{i}^{1}(t)\geq g_{i}^{2}(t)$ for all $i\in D$ and $t\geq\tau$.
Then, $u_{i}\left(  t\right)  >v_{i}\left(  t\right)  $, for all $i\in D$ and
$t\geq\tau$, where $u\left(  \text{\textperiodcentered}\right)  $,\ $v\left(
\text{\textperiodcentered}\right)  $ are the unique solutions to problem
(\ref{DelayStoch}) corresponding to $\{\phi^{1},g^{1}\}$ and $\{\phi^{2}%
,g^{2}\}$, respectively.
\end{corollary}

\section{The non-linear case}

Let us consider now the system%
\begin{align}
\frac{d}{dt}u_{i}\left(  t\right)   &  =\left[  J(t,u_{t})\right]  _{i}%
-u^{i}(t)+bu_{i}(t)(1-u_{i}(t))\frac{dw_{i}}{dt}\text{, }i\in D\text{, }%
t>\tau,\label{Nonlinear}\\
u^{i}\left(  \tau+s\right)   &  \equiv\phi^{i}\left(  s\right)  \text{, }i\in
D\text{, }s\in\lbrack-h,0]\text{.}\nonumber
\end{align}

We are now interested in proving that the components of the solution remain in
the interval $\left(  0,1\right)  $ for every moment of time. In this way, we
guarantee that the variables are probabilities if the initial conditions are
as well.

\begin{lemma}
\label{ExistSol2}Assume that $g_{r}(t)\in\lbrack0,1],$ for all $r$ and $t,$
and that $\overline{\alpha}_{i}\leq1$, for all $i$. Then, for any $\phi\in
C([-h,0],\mathbb{R}_{+}^{m})$ such that $\phi_{i}(s)\in\left(  0,1\right)  $,
for any $i\in D$ and $s\in\lbrack-h,0]$, the unique solution $u\left(
\text{\textperiodcentered}\right)  $ to (\ref{Nonlinear}) satisfies almost
surely that $u_{i}\left(  t\right)  \in\left(  0,1\right)  $ for all $i\in D$
and $t\geq\tau.$
\end{lemma}

\begin{proof}
The existence and uniqueness of local solution is again guaranteed by
\cite[Theorem 2.8, P. 154]{MaoBook}. Now we prove the statement of the Lemma.

Let $k_{0}>0$ be such that%
\[
\frac{1}{k_{0}}<\min_{s\in\lbrack-h,0]}\left\vert \phi_{i}(s)\right\vert
\leq\max_{s\in\lbrack-h,0]}\left\vert \phi_{i}(s)\right\vert <k_{0}\text{ for
all }i\in D.
\]
We define now the stopping time%
\[
\tau_{k}=\inf\{t\in\lbrack0,\infty):u_{i}\left(  t\right)  \not \in (\frac
{1}{k},1-\frac{1}{k})\text{ for some }i\in D\}.
\]
Since this sequence is increasing as $k\nearrow+\infty$, if $\tau_{\infty
}=\lim_{k\rightarrow+\infty}\tau_{k}=+\infty$ a.s., then $0<u_{i}\left(
t\right)  <1$ almost sure for $t\geq0$.

By contradiction, assume the existence of $T,\varepsilon>0$ such that
\[
P(\tau_{\infty}\leq T)>\varepsilon.
\]
In such a case there would exist $k_{1}\geq k_{0}$ such that%
\[
P(\tau_{k}\leq T)\geq\varepsilon\text{ for }k\geq k_{1}.
\]

We denote
\[
K_{0}=\{u=\left(  u_{m_{1}},...,u_{m_{2}}\right)  \in\mathbb{R}_{+}%
^{m}:0<u_{i}<1\}
\]
and define the $C^{2}$-function $V:K_{0}\rightarrow\mathbb{R}_{+}^{1}$ given
by%
\[
V\left(  u\right)  =-\sum_{i\in D}\left(  \log(1-u_{i})+\log(u_{i}\right)  ).
\]
For $\tau\leq t\leq\tau_{k}\wedge T$ we have $u\left(  t\right)  \in K_{0}$,
and then by It\^{o}'s formula we have%
\begin{align*}
&  dV(u(t))\\
&  =\sum_{i\in D}\left(  \frac{1}{1-u_{i}\left(  t\right)  }-\frac{1}{u_{i}%
}\right)  \left(  \sum_{r\in D}\int_{-h}^{0}j_{i-r}u_{r}(t+s)\alpha_{i}%
(s)d\mu(s)+\sum_{r\in\mathbb{Z}\backslash D}\int_{-h}^{0}j_{i-r}g_{r}\left(
t+s\right)  \alpha_{i}(s)d\mu(s)-u_{i}(t)\right)  dt\\
&  +\sum_{i\in D}\frac{1}{2}\left(  u_{i}^{2}\left(  t\right)  +\left(
1-u_{i}(t)\right)  ^{2}\right)  b^{2}dt+\sum_{i\in D}b(2u_{i}(t)-1)dw_{i}(t)\\
&  =I(t)dt+\sum_{i\in D}I_{i}(t)dw_{i}(t),
\end{align*}
where $I_{i}(t)=b(2u_{i}(t)-1)$. First,%
\[
\sum_{i\in D}\left(  -\frac{1}{u_{i}}\right)  \left(  \sum_{r\in D}%
j_{i-r}u_{r}\left(  t\right)  -u_{i}(t)+\sum_{r\in\mathbb{Z}\backslash
D}j_{i-r}g_{r}\left(  t\right)  \right)  \leq\sum_{i\in D}1.
\]
Second, by $\overline{\alpha}_{i}\leq1,\ 0\leq g_{r}\left(  t\right)  \leq1$
and $\left(  H2\right)  $ we obtain that%
\begin{align*}
&  \sum_{r\in D}\int_{-h}^{0}j_{i-r}u^{r}(t+s)\alpha_{i}(s)d\mu(s)+\sum
_{r\in\mathbb{Z}\backslash D}\int_{-h}^{0}j_{i-r}g_{r}\left(  t+s\right)
\alpha_{i}(s)d\mu(s)-u_{i}(t)\\
&  \leq\sum_{r\in\mathbb{Z}}j_{i-r}\int_{-h}^{0}\alpha_{i}(s)d\mu
(s)-u_{i}(t)\\
&  \leq1-u_{i}(t).
\end{align*}
Then the term $I(t)$ is estimated as follows:
\[
I\left(  t\right)  \leq\sum_{i\in D}\left(  2+b^{2}\right)  =m\left(
2+b^{2}\right)  .
\]
Integrating over $\left(  0,\tau_{k}\wedge T\right)  $ and taking expectations
we deduce that%
\begin{align*}
0  &  \leq\mathbb{E}V(u(\tau_{k}\wedge T))\leq V(\phi(0))+\mathbb{E}\int%
_{0}^{\tau_{k}\wedge T}m\left(  2+b^{2}\right)  dt+\mathbb{E}\int_{0}%
^{\tau_{k}\wedge T}\sum_{i\in D}b(2u_{i}(t)-1)dw_{i}(t)\\
&  =V(\phi(0))+m\left(  2+b^{2}\right)  \mathbb{E(}\tau_{k}\wedge T)\leq
V(\phi(0))+m\left(  2+b^{2}\right)  T.
\end{align*}

Let $\Omega_{k}=\{\omega:\tau_{k}\leq T\}$, which satisfies $P(\Omega_{k}%
)\geq\varepsilon$ for $k\geq k_{1}$. For any $\omega\in\Omega_{k}$ there
exists $i\in D$ such that either $u_{i}\left(  \tau_{k},\omega\right)
=\frac{1}{k}$ or $u_{i}\left(  \tau_{k},\omega\right)  =1-\frac{1}{k}$, so
that%
\[
V(u(\tau_{k}\wedge T,\omega))\geq\log(k)-\log(1-\frac{1}{k}).
\]
Thus,%
\[
V(\phi(0))+m\left(  2+b^{2}\right)  T\geq\mathbb{E}(1_{\Omega_{k}}V(u(\tau
_{k}\wedge T)))\geq\varepsilon\left(  \log(k)-\log(1-\frac{1}{k})\right)  .
\]
Passing to the limit as $k\rightarrow+\infty$ we arrive at a contradiction.
\end{proof}

\bigskip

\section{Asymptotic behaviour}

If we consider model (\ref{DelayStoch}) in the deterministic and autonomous
cases, that is, $b=0$ and $g_{r}(t)\equiv g_{r}\in\mathbb{R}$, and assume that
$\alpha_{i}(s)=\alpha(s)$, for all $i\in D$, then it is well known
\cite{MoVa01} that there exists a unique fixed point $\overline{u}$ given by
the solution of the system%
\begin{equation}
-M_{1}\sum_{r\in D}j_{i-r}u_{r}+u_{i}=M_{1}\sum_{r\in\mathbb{Z}\backslash
D}j_{i-r}g_{r}=b_{i},\ i\in D, \label{Fixed}%
\end{equation}
where $M_{1}=M_{1}(h):=\int_{-h}^{0}\alpha(s)d\mu(s)$, provided that
\begin{equation}
M_{1}\sum_{r\in D}j_{i-r}<1\text{ }\forall i\in D. \label{CondFixed}%
\end{equation}
Moreover, $\overline{u}_{r}\geq0$ for any $r\in D$ (see Remark 3.1 in
\cite{CarMoVa02}).

We will show that the solutions of the stochastic system remain close to this
fixed point for large times in a suitable sense.

We start with the linear case.

\begin{theorem}
\label{AsympBeh}Assume that $g_{r}\geq0$, for all $r\in D$, $b<1$ and that%
\begin{equation}
M_{1}(h)\sum_{r\in D}j_{i-r}<1-b^{2}\text{ }\forall i\in D, \label{CondM1}%
\end{equation}%
\begin{equation}
h<\frac{1}{2(1-b^{2})}\log\left(  \frac{1-b^{2}}{1-\delta(h)-b^{2}}\right)  ,
\label{Condh}%
\end{equation}
where%
\[
\delta(h)=\min_{i\in D}\left\{  1-b^{2}-M_{1}(h)\sum_{r\in D}j_{i-r}.\right\}
\]
Then, for any $\phi\in C([-h,0],\mathbb{R}_{+}^{m})$, the unique solution to
problem (\ref{Linear}) satisfies%
\[
\lim\sup_{t\rightarrow+\infty}\ \frac{1}{t}\int_{0}^{t}\sup_{\theta\in
\lbrack-h,0]}\mathbb{E}(\left\Vert u(s+\theta)-\overline{u}\right\Vert
_{\mathbb{R}^{m}}^{2}ds\leq\frac{2b^{2}\left\Vert \overline{u}\right\Vert
_{\mathbb{R}^{m}}^{2}}{\lambda^{\ast}-L(h,\lambda^{\ast})},
\]
where $L=L(h,\lambda):=2(1-\delta(h)-b^{2})e^{\lambda h}$ and $\lambda^{\ast
}\in\left(  0,2(1-b^{2})\right)  $ is such that $\lambda^{\ast}>L(h,\lambda
^{\ast}).$
\end{theorem}

\begin{remark}
It is easy to see that (\ref{Condh}) is satisfied for $h$ small enough.
Indeed, let $h_{0}$ be such that (\ref{CondM1}) holds for all $h\leq h_{0}$.
If $h<h_{0}$ is small enough such that $h<\frac{1}{2(1-b^{2})}\log\left(
\frac{1-b^{2}}{1-\delta(h_{0})-b^{2}}\right)  $, the fact that $\delta(h)$ is
non-increasing implies that $h<\frac{1}{2(1-b^{2})}\log\left(  \frac{1-b^{2}%
}{1-\delta(h)-b^{2}}\right)  .$
\end{remark}

\begin{remark}
Let $f(\lambda)=2(1-\delta(h)-b^{2})e^{\lambda h}-\lambda$. Condition
(\ref{Condh}) implies that $f(2(1-b^{2}))<0$. Then choosing $\lambda^{\ast}$
close enough to $2(1-b^{2})$ we have that $f(\lambda^{\ast})<0$, so
$\lambda^{\ast}>L(h,\lambda^{\ast}).$
\end{remark}

\begin{proof}
We define the $C^{2}$ function $V:\mathbb{R}^{m}\rightarrow\mathbb{R}$ by
\[
V(u)=\sum_{i\in D}(u_{i}-\overline{u}_{i})^{2}=\left\Vert u-\overline
{u}\right\Vert _{\mathbb{R}^{m}}^{2}.
\]
For $\lambda\in\left(  0,2(1-b^{2})\right)  $ we have%
\[
d(e^{\lambda t}V(u(t)))=\lambda e^{\lambda t}V(u(t))dt+e^{\lambda t}dV(u(t))
\]
and using Ito's formula we obtain that%
\begin{align*}
dV(u(t))  &  =\sum_{i\in D}2(u_{i}-\overline{u}_{i})\left(  \sum_{r\in D}%
\int_{-h}^{0}j_{i-r}u_{r}(t+s)\alpha(s)d\mu(s)+\sum_{r\in\mathbb{Z}\backslash
D}\int_{-h}^{0}j_{i-r}g_{r}\alpha(s)d\mu(s)-u_{i}(t)\right)  dt\\
&  +\sum_{i\in D}u_{i}^{2}(t)b^{2}dt+\sum_{i\in D}2b(u_{i}(t)-\overline{u}%
_{i})u_{i}(t)dw_{i}(t)\\
&  =I(t)dt+\sum_{i\in D}I_{i}(t)dw_{i}(t).
\end{align*}
Using (\ref{Fixed}) the term $I(t)$ is estimated by%
\begin{align*}
I(t)  &  =\sum_{i\in D}2(u_{i}-\overline{u}_{i})\left(  \sum_{r\in D}\int%
_{-h}^{0}j_{i-r}u_{r}(t+s)\alpha(s)d\mu(s)+M_{1}\sum_{r\in\mathbb{Z}\backslash
D}j_{i-r}g_{r}-u_{i}(t)\right)  +\sum_{i\in D}u_{i}^{2}(t)b^{2}\\
&  \leq\sum_{i\in D}2(u_{i}-\overline{u}_{i})\left(  \sum_{r\in D}\int%
_{-h}^{0}j_{i-r}(u_{r}(t+s)-\overline{u}_{r})\alpha(s)d\mu(s)-\left(
u_{i}(t)-\overline{u}_{i}\right)  \right)  +2b^{2}\left\Vert u-\overline
{u}\right\Vert _{\mathbb{R}^{m}}^{2}+2b^{2}\left\Vert \overline{u}\right\Vert
_{\mathbb{R}^{m}}^{2}\\
&  =\sum_{i\in D}2(u_{i}-\overline{u}_{i})\sum_{r\in D}\int_{-h}^{0}%
j_{i-r}(u_{r}(t+s)-\overline{u}_{r})\alpha(s)d\mu(s)-2(1-b^{2})\left\Vert
u-\overline{u}\right\Vert _{\mathbb{R}^{m}}^{2}+2b^{2}\left\Vert \overline
{u}\right\Vert _{\mathbb{R}^{m}}^{2}\\
&  =I_{1}(t)+I_{2}(t).
\end{align*}
Using the definition of $\delta(h)$ the first term is estimated by%
\begin{align*}
I_{1}(t)  &  \leq M_{1}(h)\sum_{i\in D}(u_{i}(t)-\overline{u}_{i})^{2}%
\sum_{r\in D}j_{i-r}+\sum_{i\in D}\sum_{r\in D}\int_{-h}^{0}j_{i-r}%
(u_{r}(t+s)-\overline{u}_{r})^{2}\alpha(s)d\mu(s)\\
&  \leq(1-\delta(h)-b^{2})\sum_{i\in D}(u_{i}(t)-\overline{u}_{i})^{2}%
+\frac{(1-\delta(h)-b^{2})}{M_{1}(h)}\int_{-h}^{0}\sum_{r\in D}(u_{r}%
(t+s)-\overline{u}_{r})^{2}\alpha(s)d\mu(s).
\end{align*}
Hence,
\begin{align*}
d(e^{\lambda t}V(u(t)))  &  \leq\lambda e^{\lambda t}\left\Vert u(t)-\overline
{u}\right\Vert _{\mathbb{R}^{m}}^{2}-2(1-b^{2})\left\Vert u-\overline
{u}\right\Vert _{\mathbb{R}^{m}}^{2}e^{\lambda t}+2b^{2}\left\Vert
\overline{u}\right\Vert _{\mathbb{R}^{m}}^{2}e^{\lambda t}+(1-\delta
(h)-b^{2})\left\Vert u(t)-\overline{u}\right\Vert _{\mathbb{R}^{m}}%
^{2}e^{\lambda t}\\
&  +\frac{(1-\delta(h)-b^{2})}{M_{1}(h)}e^{\lambda t}\int_{-h}^{0}\sum_{r\in
D}(u_{r}(t+s)-\overline{u}_{r})^{2}\alpha(s)d\mu(s)+e^{\lambda t}\sum_{i\in
D}I_{i}(t)dw_{i}(t).
\end{align*}
Integrating over $\left(  0,t\right)  $ and taking into account that
$\lambda<2(1-b^{2})$, we can discard the first two terms appearing on the
right-hand side and deduce%
\begin{align*}
e^{\lambda t}V(u(t))  &  \leq V(u(0))+\frac{2b^{2}\left\Vert \overline
{u}\right\Vert _{\mathbb{R}^{m}}^{2}}{\lambda}e^{\lambda t}+(1-\delta
(h)-b^{2})\int_{0}^{t}\left\Vert u(s)-\overline{u}\right\Vert _{\mathbb{R}%
^{m}}^{2}e^{\lambda s}ds\\
&  +\frac{(1-\delta(h)-b^{2})}{M_{1}(h)}\int_{0}^{t}e^{\lambda l}\int_{-h}%
^{0}\sum_{r\in D}(u_{r}(l+s)-\overline{u}_{r})^{2}\alpha(s)d\mu(s)dl+\int%
_{0}^{t}e^{\lambda s}\sum_{i\in D}I_{i}(s)dw_{i}(s).
\end{align*}
Taking expectations we obtain%
\begin{align*}
e^{\lambda t}\mathbb{E}\left\Vert u(t)-\overline{u}\right\Vert _{\mathbb{R}%
^{m}}^{2}  &  \leq\mathbb{E}\left\Vert u(0)-\overline{u}\right\Vert
_{\mathbb{R}^{m}}^{2}+\frac{2b^{2}\left\Vert \overline{u}\right\Vert
_{\mathbb{R}^{m}}^{2}}{\lambda}e^{\lambda t}+(1-\delta(h)-b^{2})\int_{0}%
^{t}\mathbb{E}\left\Vert u(s)-\overline{u}\right\Vert _{\mathbb{R}^{m}}%
^{2}e^{\lambda s}ds\\
&  +\frac{(1-\delta(h)-b^{2})}{M_{1}(h)}\int_{0}^{t}e^{\lambda l}\int_{-h}%
^{0}\mathbb{E}\left\Vert u(l+s)-\overline{u}\right\Vert _{\mathbb{R}^{m}}%
^{2}\alpha(s)d\mu(s)dl\\
&  \leq\mathbb{E}\left\Vert u(0)-\overline{u}\right\Vert _{\mathbb{R}^{m}}%
^{2}+\frac{2b^{2}\left\Vert \overline{u}\right\Vert _{\mathbb{R}^{m}}^{2}%
}{\lambda}e^{\lambda t}+(1-\delta(h)-b^{2})\int_{0}^{t}\mathbb{E}\left\Vert
u(s)-\overline{u}\right\Vert _{\mathbb{R}^{m}}^{2}e^{\lambda s}ds\\
&  +(1-\delta(h)-b^{2})\int_{0}^{t}\sup_{s\in\lbrack-h,0]}\mathbb{E}\left\Vert
u(l+s)-\overline{u}\right\Vert _{\mathbb{R}^{m}}^{2}e^{\lambda l}dl\\
&  \leq\mathbb{E}\left\Vert u(0)-\overline{u}\right\Vert _{\mathbb{R}^{m}}%
^{2}+\frac{2b^{2}\left\Vert \overline{u}\right\Vert _{\mathbb{R}^{m}}^{2}%
}{\lambda}e^{\lambda t}+2(1-\delta(h)-b^{2})\int_{0}^{t}\sup_{s\in
\lbrack-h,0]}\mathbb{E}\left\Vert u(l+s)-\overline{u}\right\Vert
_{\mathbb{R}^{m}}^{2}e^{\lambda l}dl.
\end{align*}
Notice that the last term makes sense thanks to \cite[Lemma 2.3, P.
150]{MaoBook}, since this implies that $u\in L^{2}(\Omega;C([-h,T];\mathbb{R}%
^{m}))\subset C([-h,T],L^{2}(\Omega;\mathbb{R}^{m}))$.

We next replace $t$ by $t+\theta$, $\theta\in\lbrack-h,0]$, $t+\theta\geq0,$
in the above inequality. Then%
\begin{align*}
&  \mathbb{E}\left\Vert u(t+\theta)-\overline{u}\right\Vert _{\mathbb{R}^{m}%
}^{2}\\
&  \leq e^{-\lambda(t+\theta)}\mathbb{E}\left\Vert u(0)-\overline
{u}\right\Vert _{\mathbb{R}^{m}}^{2}+\frac{2b^{2}\left\Vert \overline
{u}\right\Vert _{\mathbb{R}^{m}}^{2}}{\lambda}+2(1-\delta(h)-b^{2}%
)e^{-\lambda(t+\theta)}\int_{0}^{t+\theta}\sup_{s\in\lbrack-h,0]}%
\mathbb{E}\left\Vert u(l+s)-\overline{u}\right\Vert _{\mathbb{R}^{m}}%
^{2}e^{\lambda l}dl\\
&  \leq e^{-\lambda t}e^{\lambda h}\mathbb{E}\left\Vert u(0)-\overline
{u}\right\Vert _{\mathbb{R}^{m}}^{2}+\frac{2b^{2}\left\Vert \overline
{u}\right\Vert _{\mathbb{R}^{m}}^{2}}{\lambda}+2(1-\delta(h)-b^{2})e^{-\lambda
t}e^{\lambda h}\int_{0}^{t}\sup_{s\in\lbrack-h,0]}\mathbb{E}\left\Vert
u(l+s)-\overline{u}\right\Vert _{\mathbb{R}^{m}}^{2}e^{\lambda l}dl.
\end{align*}
For $t+\theta<0$ we have%
\[
e^{\lambda t}\mathbb{E}\left\Vert u(t+\theta)-\overline{u}\right\Vert
_{\mathbb{R}^{m}}^{2}\leq e^{\lambda t}\sup_{\theta\in\lbrack-h,0]}%
\mathbb{E}\left\Vert u(\theta)-\overline{u}\right\Vert _{\mathbb{R}^{m}}%
^{2}\leq e^{\lambda h}\sup_{\theta\in\lbrack-h,0]}\mathbb{E}\left\Vert
u(\theta)-\overline{u}\right\Vert _{\mathbb{R}^{m}}^{2}.
\]
Thus, if we define the function%
\[
y(t)=\sup_{\theta\in\lbrack-h,0]}\mathbb{E}\left\Vert u(t+\theta)-\overline
{u}\right\Vert _{\mathbb{R}^{m}}^{2},
\]
then
\[
e^{\lambda t}y(t)\leq y(0)e^{\lambda h}+\frac{2b^{2}\left\Vert \overline
{u}\right\Vert _{\mathbb{R}^{m}}^{2}}{\lambda}e^{\lambda t}+L(h)\int_{0}%
^{t}y(l)e^{\lambda l}dl
\]
and Gronwall's lemma implies%
\begin{align*}
e^{\lambda t}y(t)  &  \leq y(0)e^{\lambda h}+\frac{2b^{2}\left\Vert
\overline{u}\right\Vert _{\mathbb{R}^{m}}^{2}}{\lambda}e^{\lambda t}%
+L(h)\int_{0}^{t}\left(  y(0)e^{\lambda h}+\frac{2b^{2}\left\Vert \overline
{u}\right\Vert _{\mathbb{R}^{m}}^{2}}{\lambda}e^{\lambda h}e^{\lambda
l}\right)  e^{L(h)(t-l)}dl\\
&  \leq y(0)e^{\lambda h}+y(0)e^{\lambda h}(e^{L(h)t}-1)+\frac{2b^{2}%
\left\Vert \overline{u}\right\Vert _{\mathbb{R}^{m}}^{2}}{\lambda}e^{\lambda
t}+\frac{2b^{2}\left\Vert \overline{u}\right\Vert _{\mathbb{R}^{m}}^{2}%
}{\lambda}e^{\lambda h}\frac{L(h)}{\lambda-L(h)}e^{L(h)t}(e^{(\lambda
-L(h))t}-1)\\
&  \leq y(0)e^{\lambda h}e^{L(h)t}+\frac{2b^{2}\left\Vert \overline
{u}\right\Vert _{\mathbb{R}^{m}}^{2}}{\lambda-L(h)}e^{\lambda t},
\end{align*}
and, consequently,%
\[
y(t)\leq y(0)e^{\lambda h}e^{(L(h)-\lambda)t}+\frac{2b^{2}\left\Vert
\overline{u}\right\Vert _{\mathbb{R}^{m}}^{2}}{\lambda-L(h)}.
\]
Integrating over $\left(  0,t\right)  $ we have%
\begin{align*}
\frac{1}{t}\int_{0}^{t}y(s)ds  &  \leq\frac{y(0)e^{\lambda h}}{\lambda
-L(h)}\frac{1-e^{(L(h)-\lambda)t}}{t}+\frac{2b^{2}\left\Vert \overline
{u}\right\Vert _{\mathbb{R}^{m}}^{2}}{\lambda-L(h)}\\
&  \leq\frac{y(0)e^{\lambda h}}{\lambda-L(h)}\frac{1}{t}+\frac{2b^{2}%
\left\Vert \overline{u}\right\Vert _{\mathbb{R}^{m}}^{2}}{\lambda-L(h)}.
\end{align*}
Therefore,%
\[
\lim\sup_{t\rightarrow+\infty}\ \frac{1}{t}\int_{0}^{t}y(s)ds\leq\frac
{2b^{2}\left\Vert \overline{u}\right\Vert _{\mathbb{R}^{m}}^{2}}{\lambda
-L(h)},
\]
which is true when $\lambda>L(h)$. Thus, the statement is proved.
\end{proof}

\bigskip

Let us consider now system (\ref{Nonlinear}).

\begin{theorem}
Assume that $g_{r}\in\lbrack0,1]$, for all $r\in D$, $b<1$ and that
(\ref{CondM1}), (\ref{Condh}) hold. Then for any $\phi\in C([-h,0],\mathbb{R}%
_{+}^{m})$ such that $\phi_{i}(s)\in(0,1)$, for all $i\in D$ and $s\in
\lbrack-h,0]$, the unique solution to problem (\ref{Nonlinear}) satisfies%
\[
\lim\sup_{t\rightarrow+\infty}\ \frac{1}{t}\int_{0}^{t}\sup_{\theta\in
\lbrack-h,0]}\mathbb{E}(\left\Vert u(s+\theta)-\overline{u}\right\Vert
_{\mathbb{R}^{m}}^{2}ds\leq\frac{2b^{2}\left\Vert \overline{u}\right\Vert
_{\mathbb{R}^{m}}^{2}}{\lambda^{\ast}-L(h,\lambda^{\ast})},
\]
where $L=L(h,\lambda):=2(1-\delta(h)-b^{2})e^{\lambda h}$ and $\lambda^{\ast
}\in\left(  0,2(1-b^{2})\right)  $ is such that $\lambda^{\ast}>L(h,\lambda
^{\ast}).$
\end{theorem}

\begin{proof}
As before, we define the $C^{2}$-function $V:\mathbb{R}^{m}\rightarrow
\mathbb{R}$ given by
\[
V(u)=\sum_{i\in D}(u_{i}-\overline{u}_{i})^{2}=\left\Vert u-\overline
{u}\right\Vert _{\mathbb{R}^{m}}^{2}.
\]
For $\lambda\in\left(  0,2(1-b^{2})\right)  $ we have%
\[
d(e^{\lambda t}V(u(t)))=\lambda e^{\lambda t}V(u(t))+e^{\lambda t}dV(u(t)).
\]
Then, Ito's formula yields%
\begin{align*}
dV(u(t))  &  =\sum_{i\in D}2(u_{i}-\overline{u}_{i})\left(  \sum_{r\in D}%
\int_{-h}^{0}j_{i-r}u_{r}(t+s)\alpha(s)d\mu(s)+\sum_{r\in\mathbb{Z}\backslash
D}\int_{-h}^{0}j_{i-r}g_{r}\alpha(s)d\mu(s)-u_{i}(t)\right) \\
&  +\sum_{i\in D}u_{i}^{2}(t)(1-u_{i}(t)^{2}b^{2}dt+\sum_{i\in D}%
2b(u_{i}(t)-\overline{u}_{i})u_{i}(t)(1-u_{i}(t))dw_{i}(t)\\
&  =I(t)dt+\sum_{i\in D}I_{i}(t)dw_{i}(t).
\end{align*}
Since $\sum_{i\in D}u_{i}^{2}(t)(1-u_{i}(t)^{2}b^{2}\leq\sum_{i\in D}u_{i}%
^{2}(t)b^{2}$, repeating the same steps in the proof of Theorem \ref{AsympBeh}
we derive%
\begin{align*}
e^{\lambda t}V(u(t))  &  \leq V(u(0))+\frac{2b^{2}\left\Vert \overline
{u}\right\Vert _{\mathbb{R}^{m}}^{2}}{\lambda}e^{\lambda t}+(1-\delta
(h)-b^{2})\int_{0}^{t}\left\Vert u(s)-\overline{u}\right\Vert _{\mathbb{R}%
^{m}}^{2}e^{\lambda s}ds\\
&  +\frac{(1-\delta(h)-b^{2})}{M_{1}(h)}\int_{0}^{t}e^{\lambda l}\int_{-h}%
^{0}\sum_{r\in D}(u_{r}(l+s)-\overline{u}_{r})^{2}\alpha(s)d\mu(s)dl+\int%
_{0}^{t}e^{\lambda s}\sum_{i\in D}I_{i}(s)dw_{i}(s)
\end{align*}
and, taking expectations, we infer
\[
e^{\lambda t}\mathbb{E}\left\Vert u(t)-\overline{u}\right\Vert _{\mathbb{R}%
^{m}}^{2}\leq\mathbb{E}\left\Vert u(0)-\overline{u}\right\Vert _{\mathbb{R}%
^{m}}^{2}+\frac{2b^{2}\left\Vert \overline{u}\right\Vert _{\mathbb{R}^{m}}%
^{2}}{\lambda}e^{\lambda t}+2(1-\delta(h)-b^{2})\int_{0}^{t}\sup_{s\in
\lbrack-h,0]}\mathbb{E}\left\Vert u(l+s)-\overline{u}\right\Vert
_{\mathbb{R}^{m}}^{2}e^{\lambda l}dl.
\]
Again, the last term in the previous equation is well defined thanks
\cite[Lemma 2.3, P. 150]{MaoBook}. Notice that as the solutions belong to
$(0,1)$ almost surely, the term in front of the noise has sub-linear growth
Lemma 2.3 in \cite{MaoBook} can be applied to this nonlinear case. The rest of
the proof repeats the same argument as in Theorem~\ref{AsympBeh}.
\end{proof}

\section{Application to Life Tables}

In this section, we will appply model (\ref{Nonlinear}) to predict the
probability of death by age in Spain.

\subsection{Life Tables: mortality modeling through a stochastic delay
approach}

Life tables are among the most widely used tools for studying survival and
mortality patterns in a population. In demography, for instance, mortality
constitutes one of the terms of the component method
(\cite{componentes-eurostat}, \cite{componentes-ine}) used for population
estimates. In actuarial science, particularly in insurance applications, life
tables are a fundamental tool for calculating, for example, adverse scenarios
that insurance companies must be prepared to face.

These tables are structured around interrelated biometric functions (see, for
example, \cite{Chiang}, \cite{BenjaminPollard2}, \cite{AyusoCoGui}). Among
these functions, we can highlight the probability of surviving to a completed
age $x$, $p_{x}$; its complement, the probability of death, $q_{x}$; life
expectancy at age $x$, $e_{x}$; and the central death rate by age, $m_{x}$.

In practice, the true values of these functions are generally unknown and must
be estimated from observed data. This estimation process can be approached
from different methodological perspectives, typically grouped into three main
categories: (i) stochastic versus non-stochastic models, (ii) parametric
versus non-parametric models, and (iii) static versus dynamic approaches. Each
of these axes defines key aspects in the modeling process: whether randomness
is considered, whether structured functional forms are imposed, and whether
the temporal evolution of mortality is incorporated.

Traditionally, the estimation of age-specific death probabilities ($q_{x}$)
has been performed using data from a single time period. This procedure
generates so-called crude death rates, which may lack desirable smoothness
properties, such as the expected continuity between adjacent ages. To correct
these irregularities, classical laws such as Gompertz's law \cite{Gompertz},
Gompertz-Makeham's law \cite{Makeham}, or the Heligman-Pollard model
\cite{Helligman}, among others, have been proposed.

However, these approaches generally adopt a parametric and static perspective,
in which the function $q_{x}^{t}$ is assumed to remain constant over nearby
time intervals. Nevertheless, it is well known that mortality rates evolve
over time, meaning that assuming $q_{x}^{t}=q_{x}^{t+h}$ can lead to
significant errors in many applications, such as pension expenditure
projections or the calculation of technical reserves for insured portfolios.
This realization has driven the development of dynamic mortality models, where
the temporal evolution of $q_{x}^{t}$ is modeled explicitly.

Among the best-known and most widely used dynamic models are the Lee-Carter
models \cite{LeeCarter}, the CBD model \cite{Cairns}, and the extended M3--M7
family of models \cite{Cairns2009}, \cite{Cairns2011}. These models introduce
temporal improvement factors and stochastic components that capture the
uncertainty associated with future mortality, offering interpretable
structures and reasonable predictive performance.

In line with these models, in \cite{MoVa14} a non-parametric dynamic model,
based on kernel smoothing techniques to estimate mortality rates, was
proposed. This model avoids rigid functional assumptions and uses a system of
non-local differential equations to approximate the evolution of $q_{x}^{t}$
over time. Although it successfully reproduces qualitative features of
observed mortality, its predictive horizon is short (two to three years) due
to the absence of historical information. To overcome this limitation,
\cite{MoVa01} proposed an improved model that incorporates past information
through a delay term, specifically using mortality improvement rates
\cite{Cairns}. This formulation retains the dynamic and non-parametric
character of the original model while substantially enhancing its predictive
capacity, extending its utility to time horizons of five to ten years. The
underlying idea is that mortality trajectories are partially path-dependent,
and thus the historical evolution must be considered.

This delay-based model remains within a non-stochastic framework, using
observed improvement rates to incorporate past dynamics. It offers a robust
alternative to stochastic models when precise deterministic prediction is
required, as in regulatory contexts (e.g., Solvency II, \cite{Solvencia 2}).
However, despite its improved predictive performance, a major limitation
persists: it does not offer 'alternative' mortality evolution scenarios, nor
does it allow the calculation of quantiles or the construction of confidence intervals.

To address this issue, the non-local model proposed in \cite{MoVa14} was
extended into a stochastic model by introducing a random term into the
original system of non-local differential equations. This extension aimed to
capture both the temporal evolution and the intrinsic variability of the
mortality phenomenon. The resulting model \cite{CarMoVa02} strikes a
compelling balance between interpretability, flexibility, and robustness, and
aligns with the family of stochastic, non-parametric, and dynamic models (see
\cite{Brouhns}, \cite{Bijak}, \cite{CopasHaber}).

Despite providing a novel approach to integrating stochastic variability and
non-parametric smoothing within a single mortality modeling framework, the
need to incorporate historical information about mortality evolution became
apparent. This leads to the model proposed in the present work, which can be
classified as a dynamic, non-parametric, and stochastic model that also
integrates historical information to capture the temporal evolution of mortality.

\subsection{Dynamical kernel graduation: combining stochastic behavior and
historical data}

The model proposed in equation (\ref{DelayStoch}) contains both delay and
stochastic terms. In a similar way as in \cite{MoVa01}, in this work the delay
term takes into account the history using the concept of \emph{Improvement
Rate} (see also \cite{DoddImprovRates}, \cite{HabermanRensahImprovRat},
\cite{HabermanRenshahImprovRat2} and \cite{PremioMapfre2015}). These rates,
denoted by $r_{x}^{t_{0},t_{1}},$ treat to characterize the evolution of the
mortality year-to-year or between two arbitrary moments, $t_{0}$ and $t_{1}$,
at a concrete age $x;$ they are defined by $r_{x}^{t_{0},t_{1}}=q_{x}^{t_{1}%
}/q_{x}^{t_{0}}$. The difference $d=t_{1}-t_{0}>0$ is the delay. Using these
improvement rates, we define the \emph{global improvement rate, }denoted by
$\overline{r}_{t_{0}}^{t_{1}}$, as the coefficient of the linear model
(without constant term) of the observed death rates, that is, we fit the
linear model
\[
q_{\cdot}^{t_{1}}=\overline{r}_{t_{0}}^{t_{1}}\cdot q_{\cdot}^{t_{0}}%
\]
with the data $\left\{  q_{0}^{t_{0}},\ldots,q_{\omega}^{t_{0}}\right\}  $ and
$\left\{  q_{0}^{t_{1}},\ldots,q_{\omega}^{t_{1}}\right\}  $.

The procedure can be summarized as in (\cite{MoVa01}):

\begin{enumerate}
\item We consider the observed mortality rates at each age ($x$) and each year
($t$): $q_{x}^{t}$, $x\in D=\left\{  0,\ldots,100\right\}  $, $t\in\left\{
1908,\ldots,2018\right\}  $. Also, we consider the values of $g_{x}\left(
t\right)  $, which is the rate of death either at "negative ages" or after the
\textit{actuarial infinite (we have chosen it equal to }$100$\textit{)}.

\item We estimate the improvement rates for each age and delay, that is:%
\[%
\begin{array}
[c]{cccccccc}%
\frac{{\small year}}{delay} & 1 & 2 & 3 & 4 & ... & 109 & 110\\
& \overline{r}_{1908}^{1909}\sim\frac{q_{\cdot}^{1909}}{q_{\cdot}^{1908}} &
\overline{r}_{1908}^{1910}\sim\frac{q_{\cdot}^{1910}}{q_{\cdot}^{1908}} &
\overline{r}_{1908}^{1911}\sim\frac{q_{\cdot}^{1911}}{q_{\cdot}^{1908}} &
\overline{r}_{1908}^{1912}\sim\frac{q_{\cdot}^{1912}}{q_{\cdot}^{1908}} &
\ldots & \overline{r}_{1908}^{2017}\sim\frac{q_{\cdot}^{2017}}{q_{\cdot
}^{1908}} & \overline{r}_{1908}^{2018}\sim\frac{q_{\cdot}^{2018}}{q_{\cdot
}^{1908}}\\
& \overline{r}_{1909}^{1910}\sim\frac{q_{\cdot}^{1910}}{q_{\cdot}^{1909}} &
\overline{r}_{1909}^{1911}\sim\frac{q_{\cdot}^{1911}}{q_{\cdot}^{1909}} &
\overline{r}_{1909}^{1912}\sim\frac{q_{\cdot}^{1912}}{q_{\cdot}^{1909}} &
\overline{r}_{1909}^{1913}\sim\frac{q_{\cdot}^{1913}}{q_{\cdot}^{1909}} &
\ldots & \overline{r}_{1909}^{2017}\sim\frac{q_{\cdot}^{2017}}{q_{\cdot
}^{1909}} & -\\
& \overline{r}_{1910}^{1911}\sim\frac{q_{\cdot}^{1911}}{q_{\cdot}^{1910}} &
\overline{r}_{1910}^{1912}\sim\frac{q_{\cdot}^{1912}}{q_{\cdot}^{1910}} &
\overline{r}_{1910}^{1913}\sim\frac{q_{\cdot}^{1913}}{q_{\cdot}^{1910}} &
\overline{r}_{1910}^{1914}\sim\frac{q_{\cdot}^{1914}}{q_{\cdot}^{1910}} &
\ldots & - & -\\
& \overline{r}_{1911}^{1912}\sim\frac{q_{\cdot}^{1912}}{q_{\cdot}^{1911}} &
\overline{r}_{1911}^{1913}\sim\frac{q_{\cdot}^{1913}}{q_{\cdot}^{1911}} &
\overline{r}_{1911}^{1914}\sim\frac{q_{\cdot}^{1914}}{q_{\cdot}^{1911}} &
\overline{r}_{1911}^{1915}\sim\frac{q_{\cdot}^{1915}}{q_{\cdot}^{1911}} &
\ldots & - & -\\
& \vdots & \vdots & \ldots & - & - & - & -\\
& \overline{r}_{2016}^{2017}\sim\frac{q_{\cdot}^{2017}}{q_{\cdot}^{2016}} &
\overline{r}_{2016}^{2018}\sim\frac{q_{\cdot}^{2018}}{q_{\cdot}^{2016}} & - &
- & - & - & -\\
& \overline{r}_{2017}^{2018}\sim\frac{q_{\cdot}^{2018}}{q_{\cdot}^{2017}} &
- & - & - & - & - & -
\end{array}
\]

\item We consider the modified spheric function $\widetilde{\gamma}$ given by%
\[
\widetilde{\gamma}\left(  s\right)  =\left\{
\begin{array}
[c]{c}%
T\text{, if }s<b-a,\\
\gamma\left(  b-s\right)  \text{, if }b-a\leq s\leq b,\\
0\text{, if }s>b,
\end{array}
\right.
\]
where
\[
\gamma\left(  s\right)  =\left\{
\begin{array}
[c]{c}%
\frac{T}{2}\left\{  3\frac{s}{a}-\left(  \frac{s}{a}\right)  ^{3}\right\}
\text{, if }0\leq s\leq a,\\
T\text{, if }s>a,
\end{array}
\right.
\]
with $a=20$, $b=30$ and $T=1$. We will use this function to calculate the
annual improvement rates by delay, $\overline{\overline{r}}^{d}$, which
summarizes the previous information of the mortality process. We calculate a
ponderated mean, with respect to the delay, of the improvement rates. For this
aim, we define the vector $v=(v_{i})$, $i=1,...,n-1$, $n=111$, given by%
\[
v_{i}=\widetilde{\gamma}\left(  i\right)  .
\]
Then, for each delay $d=1,...,110$, we define the vector $v^{d}=(v_{i}^{d})$,
$i=1,...,n-d,$ by
\[
v_{i}^{d}=\frac{v_{i}}{\sum_{j=1}^{n-d}v_{j}},
\]
and the%
\begin{equation}
\overline{\overline{r}}^{d}=\sum_{i=1}^{n-d}v_{i}^{d}\overline{r}%
_{\overline{t}+i}^{\overline{t}+i+d}, \label{global_ImprovRates}%
\end{equation}
being $\overline{t}=1907.$ In this way, the values $\overline{r}_{t_{0}%
}^{t_{1}}$ which are far enough from the present moment of time are not taken
into account. We will refer to this values as the \emph{global improvement
rates for the delay }$d$.

\item Using the improvement rates by delay, $\overline{\overline{r}}^{d}$, we
define the function $\alpha\left(  s\right)  $ from (\ref{DelayStoch}), by
putting%
\begin{equation}
\alpha\left(  -s\right)  =1+\beta s\text{, }s\in\lbrack0,h],
\label{AlfaFunction}%
\end{equation}
where $\beta$ is the coefficient of the linear regression obtained from the
data $\overline{\overline{r}}^{d}$, $d=1,2,...,h$, and $h\leq110$ is the
maximum delay to be considered.

\item The experience of studying the mortality phenomenon allows us to assure
that the importance of these rates are not the same for all delays. Indeed,
the importance of the improvement rates increases when they are close to the
time of prediction. To take this into account, we assume that the importance
of these rates is modulated by a probability distribution\ function. To do
this, we consider the exponential probability density function with the form%
\[
f_{\lambda}\left(  s\right)  =\left\{
\begin{array}
[c]{c}%
\lambda e^{-\lambda s}\text{ if }s\geq0\text{,}\\
0\text{ if }s<0,
\end{array}
\right.
\]
and obtain a density function in the interval $[0,h]$ by putting%
\begin{equation}
\widehat{f}_{\lambda}\left(  s\right)  =\frac{\lambda e^{-\lambda s}}{\int%
_{0}^{h}\lambda e^{-\lambda s}ds},\ s\in\lbrack0,h]. \label{Exp2}%
\end{equation}
Then we define the density function $\xi\left(  \text{\textperiodcentered
}\right)  $ from (\ref{DelayStoch}) by $\xi\left(  s\right)  =\widehat{f}%
_{\lambda}\left(  -s\right)  $ for $s\in\lbrack-h,0].$

In the particular case when in the numerical approximations we consider only
integer delays, we can discretize the interval $[0,h]$ by using a finite
number of integer delays $s=\left\{  0,1,\ldots,d_{\max}\right\}  $, where
$d_{\max}=h\leq n-1$, which is the maximum delay to be considered. Thus,
instead of (\ref{Exp2}) we will use the discretized probability function
\begin{equation}
f^{\ast}\left(  s\right)  =\frac{\hat{f}_{\lambda}\left(  s\right)  }%
{\sum_{s=0}^{d_{\max}}\hat{f}_{\lambda}\left(  s\right)  },
\label{probabilidad_exponencial}%
\end{equation}
which approximates (\ref{Exp2}) at $s=0,1,...,d_{\max}$.

Using the annual improvement rates, and the exponential distribution, we can
define the\textit{\ weighted improvement rates} as%
\begin{equation}
R_{weight}^{d}=\overline{\overline{r}}^{d}\cdot f^{\ast}\left(  d\right)  ,
\label{weighted_improvementRates}%
\end{equation}
but they are not used in the numerical method.
\end{enumerate}

In a similar way as in previous works (see \cite{MoVa14}, \cite{MoVa01}), for
each age $x$, for an arbitrary moment $t$ and for a time step $\tau>0$, the
probability of death at $t+\tau$, denoted by $q_{x}\left(  t+\tau\right)  $,
depends on:

\begin{itemize}
\item all graduate values at moment $t$, $q_{z}\left(  t\right)  $, $z\in
D\subset\mathbb{Z}$ (via Gaussian kernel graduation, see \cite{AyusoCoGui}).

\item all previous moments of time $t+s$, $s\in\left[  -h,0\right]  $ (via
improvement rates).
\end{itemize}

In the real world, when a phenomenon has a random nature, that is, there
exists some kind of noise which can be intrinsic to the process under study,
it is more suitable to introduce random fluctuations, for example to forecast
adverse scenarios. Then, in this work, in a similar way as in \cite{CarMoVa02}%
, we introduce the stochastic term $bq_{i}(t)(1-q_{i}(t))\dfrac{dw_{i}}{dt}$
in the model, giving rise to equation (\ref{Nonlinear}). Applying the
Euler-Maruyama discrete time approximation \cite{KloedenPlaten}, the relation
between $q_{x}\left(  t+\tau\right)  $ and $\left\{  q_{z}\left(  t\right)
\right\}  $, $\left\{  q_{z}\left(  t+s\right)  \right\}  _{s\in\left[
-h,0\right[  }\,$\ becomes:%
\begin{align}
q_{x}\left(  t+\tau\right)   &  =\frac{1}{1+\tau}q_{x}\left(  t\right)
+\frac{\tau}{1+\tau}\sum_{z\in D}\int_{-h}^{0}j_{x-z}q_{z}\left(  t+s\right)
\alpha\left(  s\right)  d\mu\left(  s\right)  \label{exp_discreta_con_retardo}%
\\
&  +\frac{\tau}{1+\tau}\sum_{z\in\mathbb{Z}\backslash D}\int_{-h}^{0}%
j_{x-z}g_{z}\left(  t+s\right)  \alpha\left(  s\right)  d\mu\left(  s\right)
+bq_{x}(t)(1-q_{x}(t))(W_{x,t+\tau}-W_{x,t}),\nonumber
\end{align}
where $j$ is a suitable kernel (in this work a Gaussian kernel), $g_{z}\left(
\cdot\right)  $ is the rate of death either at "negative ages" or after the
\textit{actuarial infinite, }$d\mu(s)=\widehat{f}_{\lambda}\left(  -s\right)
ds$, where $\widehat{f}_{\lambda}(-s),\alpha\left(  s\right)  $ are defined
above. Relation (\ref{exp_discreta_con_retardo}) is consistent with the
empirical experience on the\ actuarial practice.

In order to evaluate the integral $\int_{-h}^{0}j_{x-z}g_{z}\left(
t+s\right)  \alpha\left(  s\right)  d\mu\left(  s\right)  $ we will use the
classical Riemann sum with time step $1$ and aproximate the function
$\widehat{f}_{\lambda}(s)$ using the discretized probability function
(\ref{probabilidad_exponencial}).

We observe that when the values $q_{x}(t)$ decrease over time, which is our
case, the coefficient $\beta$ in (\ref{AlfaFunction}) is negative. Then%
\[
\overline{\alpha}=\int_{-h}^{0}\alpha\left(  s\right)  \widehat{f}_{\lambda
}\left(  s\right)  ds=\int_{-h}^{0}\left(  1-\beta s\right)  \widehat{f}%
_{\lambda}\left(  s\right)  ds=1-\beta\int_{-h}^{0}s\widehat{f}_{\lambda
}\left(  s\right)  ds\leq1,
\]
so that assumption $\overline{\alpha}\leq1$ is satisfied in Lemma
\ref{ExistSol2}.

\subsection{Numerical simulation}

We will implement our nonlinear stochastic model with delay
(\ref{exp_discreta_con_retardo}) (NLSD for short)\ to predict the probability
of death in Spain.

\subsubsection{Data and parameters}

The dataset used in this work has been obtained from \cite{HMD}. The variables
are the population and the central mortality rates for each age, which are
taken from $0$ to $100$ years old (actuarial infinity). We use the observed
values in Spain in the period $1908-2023$.

The dataset has been splitted in two subsets. First, the period $1908-2018$ is
used to fit and calibrate the model; second, the period $2019-2023$ is used
for the validation of the model.

We have chosen the value $\lambda=\frac{11}{12}$ in the function (\ref{Exp2}).

The maximum delay $h$ was set to $90$. With this value, the estimated slope
$\beta$ is $-0.003473.$

As a kernel we have chosen a discrete Gaussian Kernel with $0$ mean, variance
equal to $1$ and bandwidth $b_{w}=0.25,$ which is defined as follows. We
consider a finite set of ages $A=\left\{  -m,-m+1,...,x,\ldots,M-1,M\right\}
$, where $0\leq x\leq100,\ m\geq0,\ M\geq100$. We define the set of distances
\[
\mathbb{D}_{x}\mathbb{=}\left\{  d_{-m},d_{-m+1},\ldots,d_{0},d_{1}%
,\ldots,d_{M}\right\}  =\mathbb{\{}x+m,x+m-1,...,1,0,1,...,M-x-1,M-x\}.
\]
Then, we define the truncated gaussian kernel $\widehat{K}_{b}\left(
\text{\textperiodcentered}\right)  $ as%
\[
\widehat{K}_{b}\left(  k\right)  =\frac{K_{b}\left(  k\right)  }{\sum_{\xi
\in\mathbb{D}}K_{b}\left(  \xi\right)  }\text{,\ }k\in\mathbb{D}\text{,}%
\]
where $K\left(  \cdot\right)  $ denotes a density function from a standard
gaussian random variable and $K_{b}(x)=\frac{1}{b_{w}}K(\frac{x}{b_{w}})$,
$b_{w}=0.25$.

In expression (\ref{exp_discreta_con_retardo}) we consider the truncated
summatories with the restriction $-50\leq z\leq150$. Thus, we set
$m=50,\ M=150$ for the discrete Gaussian kernel. In this way, for
$x\in\{0,1,...,100\}$ we obtain:%
\[
j_{x-z}=\frac{K_{b}\left(  \left\vert x-z\right\vert \right)  }{\sum_{\xi
\in\mathbb{D}_{x}}K_{b}\left(  \xi\right)  },\text{ for }z\in\{-50,...,150\},
\]
where $\mathbb{D}_{x}=\{x+50,...,150-x\}.$

Related with \textit{actuarial infinity (}$100$ years old)\textit{, }the
values of $g_{x}$ for $x>100$ are taken to be equal to $0.385$, in a similar
way as in \cite{PremioMapfre2015}. For $x<0$ the values of $g_{x}^{t}$ are
estimated using the expression $\frac{3}{4}\overset{\circ}{q}_{0}^{t}+\frac
{1}{4}\overset{\circ}{q}_{1}^{t}$.

The proposed method enables forecasting over an arbitrary time horizon. Also,
the method makes it possible to obtain several trajectories, that is, an
ensemble of predictions. In particular, we have considered a time horizon of
$15$ years and a number of the trajectories equal to $500$.

The time step $\tau$ in (\ref{exp_discreta_con_retardo}) is taken equal to
$1$. Although it is possible to use a smaller value for $\tau$ by
interpolating the values of the variable in the past, the results are rather similar.

The increments $W_{x,t+\tau}-W_{x,\tau}$ are obtained using the Box-Muller
algorythm \cite{KloedenPlaten}, which gives a pair of pseudo-random numbers
$\left(  z_{1},z_{2}\right)  $ that are independent and normally distributed
with zero mean and variance equal to $\tau$.

We have obtained the predicted trajectories for several values of the
parameter $b$ determining the intensity of the noise: $b=0$.$1,0$.$05$ and
$0$.$025.$

\subsubsection{Indicators}

To validate the method and to determine if the proposed technique can be used
in real applications, we define several indicators. These indicators, also,
are used to compare different models. We consider two types of measures which
can be classified as error, count and central measures.

\textbf{Error measures. }These measures compare the observed mortality rates
with a synthetic trajectory. In this case, we use the \textit{mean (or median)
trajectory }of the realizations. Then, we calculate the \textit{mean quadratic
difference, }$I_{MqD}^{t}$, or the \textit{mean relative quadratic
difference}, $I_{MRqD}^{t}$, for each year. In particular we use the
expressions%
\begin{equation}
I_{MqD}^{t}=\frac{\sum_{x=0}^{100}\left(  q_{x}^{t,obs}-\mathbb{E}(q_{x}%
^{t})\right)  ^{2}}{101}\text{ or }\chi_{_{MRqD}}^{t}=\frac{\sum_{x=0}%
^{100}\frac{\left(  q_{x}^{t,obs}-\mathbb{E}(q_{x}^{t})\right)  ^{2}%
}{\mathbb{E}(q_{x}^{t})}}{101}\text{. } \label{error_measures}%
\end{equation}
Here, $q_{x}^{t,obs}$ are the observed rates to age $x$ at time $t$, and
$\mathbb{E}(q_{x}^{t})$ are the mean value of the trajectories of the computed realizations.

\textbf{Count measures. }In the stochastic paradigm, it can be appropiate to
use other indicators to determine if a method is good. The model proposed in
this work, as well as other models which are used in the validation step,
allows us to construct s\textit{ynthetic empirical confidence intervals}
$IC_{1-\alpha}$ with a given $\alpha$ level. Then, we define several
indicators using these confidence intervals. In particular, we put%
\begin{equation}
I_{c,1-\alpha}^{t}=\#\left\{  x:\text{ }q_{x}^{t,obs}\not \in IC_{1-\alpha
}\right\}  \text{.} \label{count_measures}%
\end{equation}
For each year $I_{c,1-\alpha}^{t}$ summarizes the number of ages of the
observed data at year $t$ that do not belong to the $\alpha-$synthetic
confidence interval, $\alpha\in\lbrack0,1]$. We will use $1-\alpha=0.98,0.90$
and $0.80$.

\textbf{Central measures and variability. }It is important to point out that a
stochastic model is suitable if it achieves a good balance between coverage
and precision. For instance, if the confidence intervals are narrow but do not
contain the observed (or future) values, such a model underestimates
uncertainty and may lead to serious consequences---for example, if an
insurance company fails to allocate sufficient capital reserves to meet future
claims. On the other hand, if the resulting confidence intervals are too wide,
even if they always include the observed or future values, the model
overestimates uncertainty, thus losing predictive value and potentially
causing significant harm---for example, by requiring capital to be reserved
for specific purposes, thereby limiting its availability for others, such as
healthcare or pensions. Hence, in the indicators of this type we take into
account both the precision of the mean values and the dispersion of the
realizations in order to compare the methods. Using the same notation as in
(\ref{error_measures}) we define, in a similar way as in \cite{DawidSebastini}%
, the following indicator:%
\begin{equation}
I_{CT,\tau}^{t}=\frac{\sum_{x=0}^{100}\left(  q_{x}^{t,obs}-\mathbb{E}%
(q_{x}^{t})\right)  ^{2}}{\left(  \sigma_{x}^{t}\right)  ^{\tau}},
\label{Central1}%
\end{equation}
where $\sigma_{x}^{t}$ is the standard deviation of the trajectories of the
computed realizations to age $x$ at time $t$ and either $\tau=1$ or $\tau=2.$

\subsubsection{Software}

The software used to implement the numerical method has been MATLAB
(versi\'{o}n R2024b).

The R-software (\cite{R-software}) has been used to download the dataset from
\cite{HMD} (using package \textit{demography, \cite{demography}}). Also, the
R-software has been used to implement the models Lee-Carter, Renshaw-Haberman,
CBD and the family M3-M7 using the package StMoMo (\cite{Villegas}); this
package has been used to determine the best models. The R-software has been
used to estimate the count, error and central measures used to validate our
model and to determine which are the best models. The figures have been
created using the R-software.

\subsubsection{Numerical results}

This section is dedicated to the validation of the proposed method and to
compare it with others techniques.

In this first part, we show graphically how the NLSD method reproduces the
qualitative bahavior of the mortality curve. In the second part, we compare
the NLSD method with the classical ones, in particular with the Lee-Carter and
Renshaw-Haberman methods.

\begin{figure}[H]
\centerline{\scalebox{0.5}{\includegraphics[angle=0]{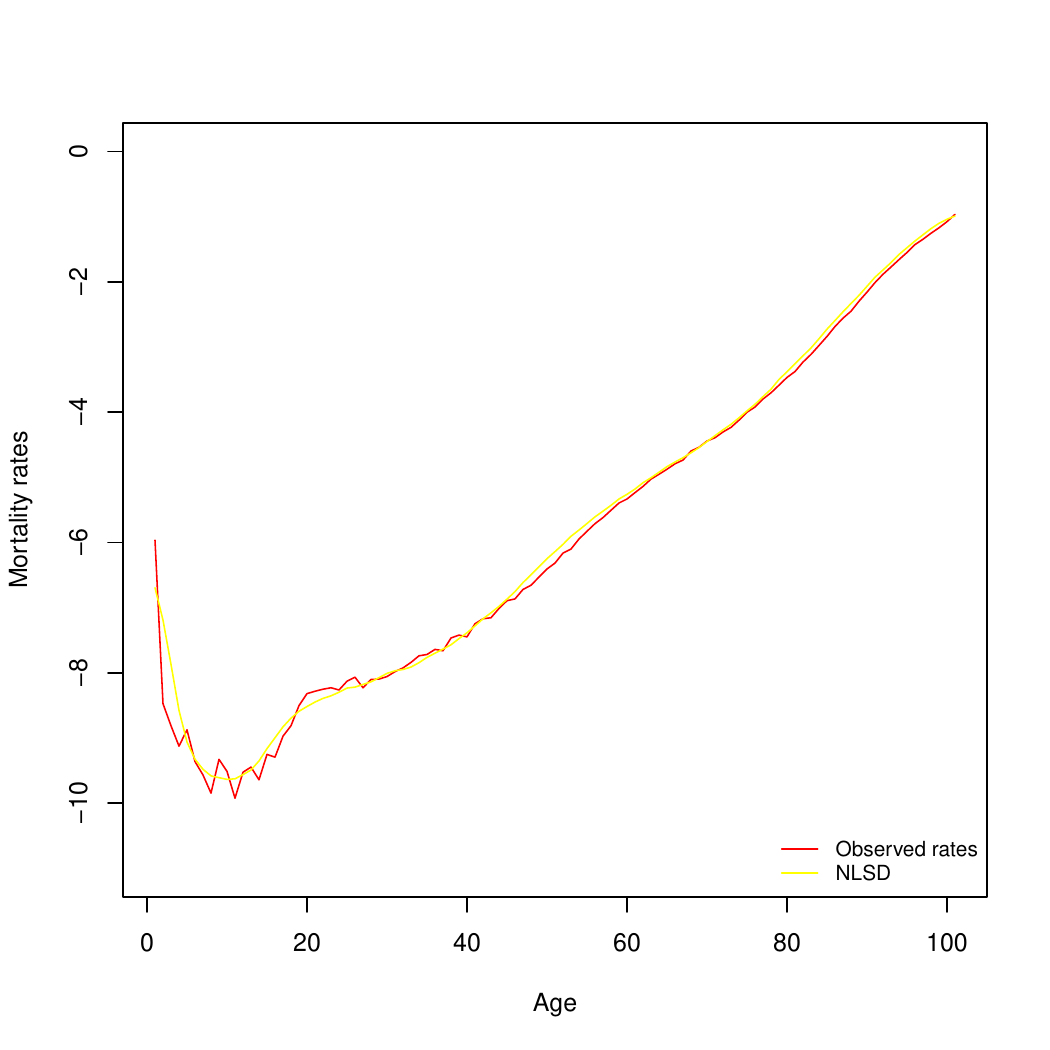}}}\caption{{\protect\footnotesize {Mean
trajectory: $b=0.1$}}}%
\label{Fig1}%
\end{figure}

Figure \ref{Fig1} shows, for the year 2023, that the mean trajectory obtained
by the NLSD method is closed to the observed rates. The same behavior can be
verified for the rest of the years in the validation period ($2019-2023$).
Also, qualitatively, we can observe how the mean realization reproduces the
form of the mortality curve, with the usual parts: adaptation to the
environment (ages 0-16), natural longevity (ages 16-100) and social jump (ages 16-27).

\begin{figure}[H]
    \centering

    \begin{subfigure}[b]{0.6\textwidth}
        \centering
        \includegraphics[scale=0.5,width=\textwidth]{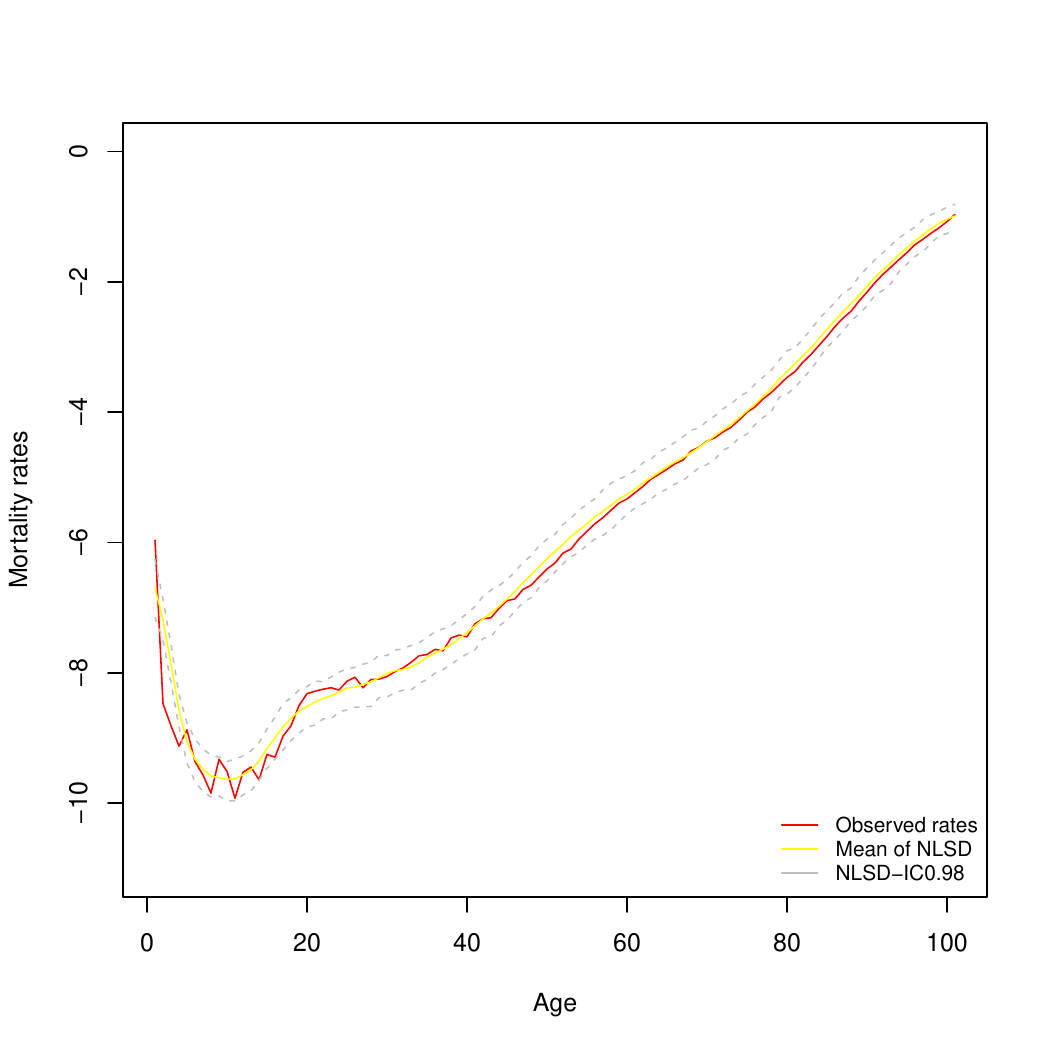}
        \caption{{Confidence Level: $1-\alpha=0.98$}}
        \label{Fig21}
    \end{subfigure}

    \vspace{0.5cm} 

    \begin{subfigure}[b]{0.45\textwidth}
        \centering
        \includegraphics[scale=0.5,width=\textwidth]{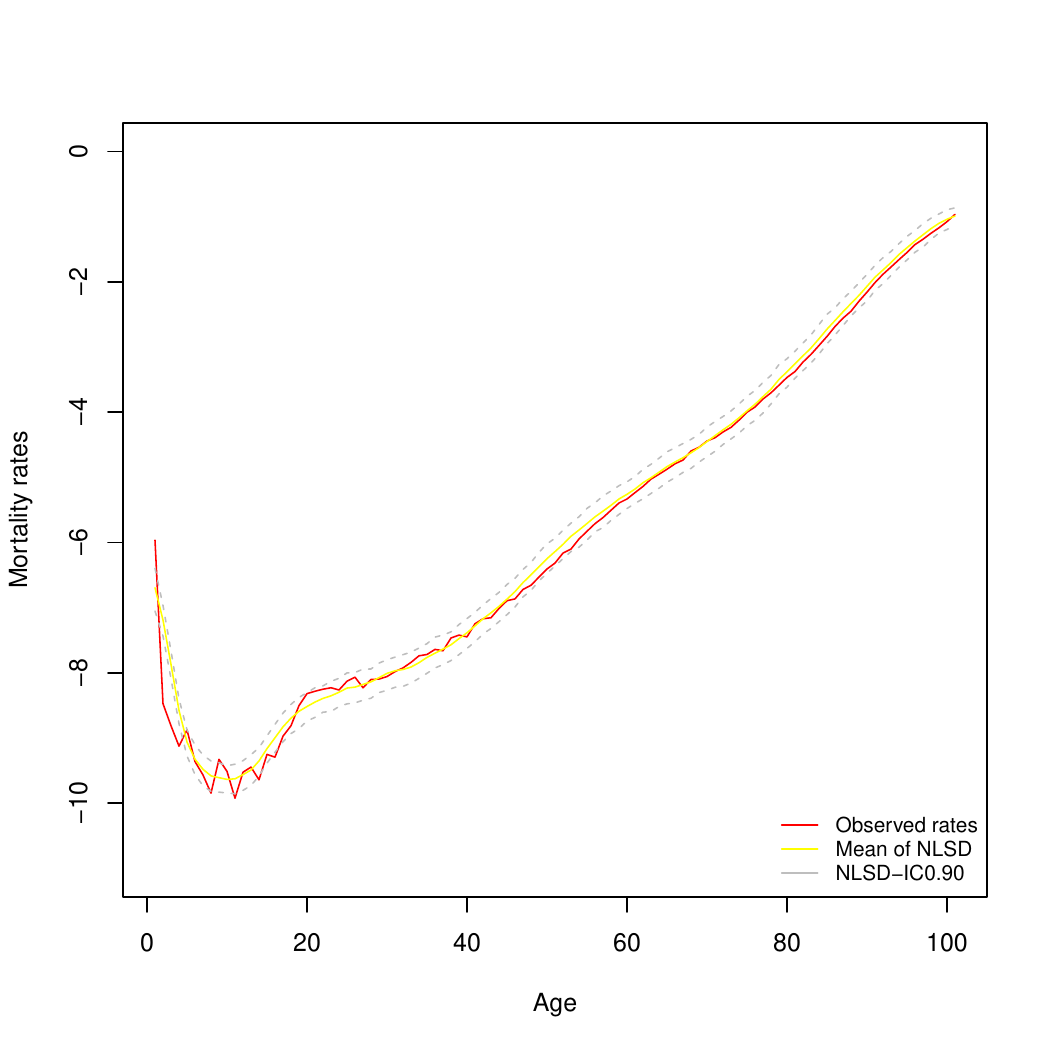}
        \caption{{Confidence Level: $1-\alpha=0.9$}}
        \label{Fig22}
    \end{subfigure}
    \hfill
    \begin{subfigure}[b]{0.45\textwidth}
        \centering
        \includegraphics[scale=0.5,width=\textwidth]{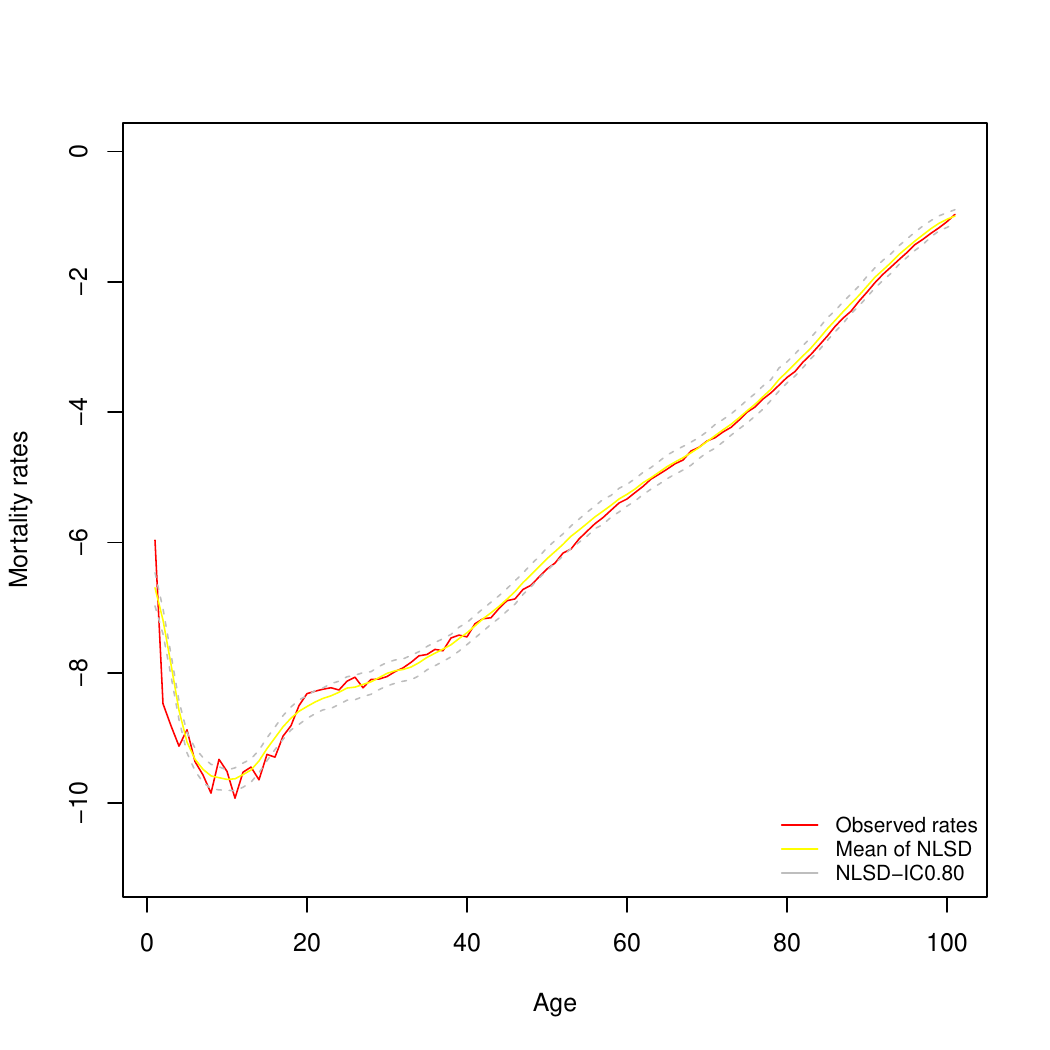}
        \caption{{Confidence Level:$1-\alpha=0.8$}}
        \label{Fig23}
    \end{subfigure}

    \caption{Confidence Intervals for several confidence levels ($b=0.1$).}
    \label{fig:figura_CIb01}
\end{figure}

\begin{figure}[H]
    \centering

    \begin{subfigure}[b]{0.6\textwidth}
        \centering
        \includegraphics[scale=0.5,width=\textwidth]{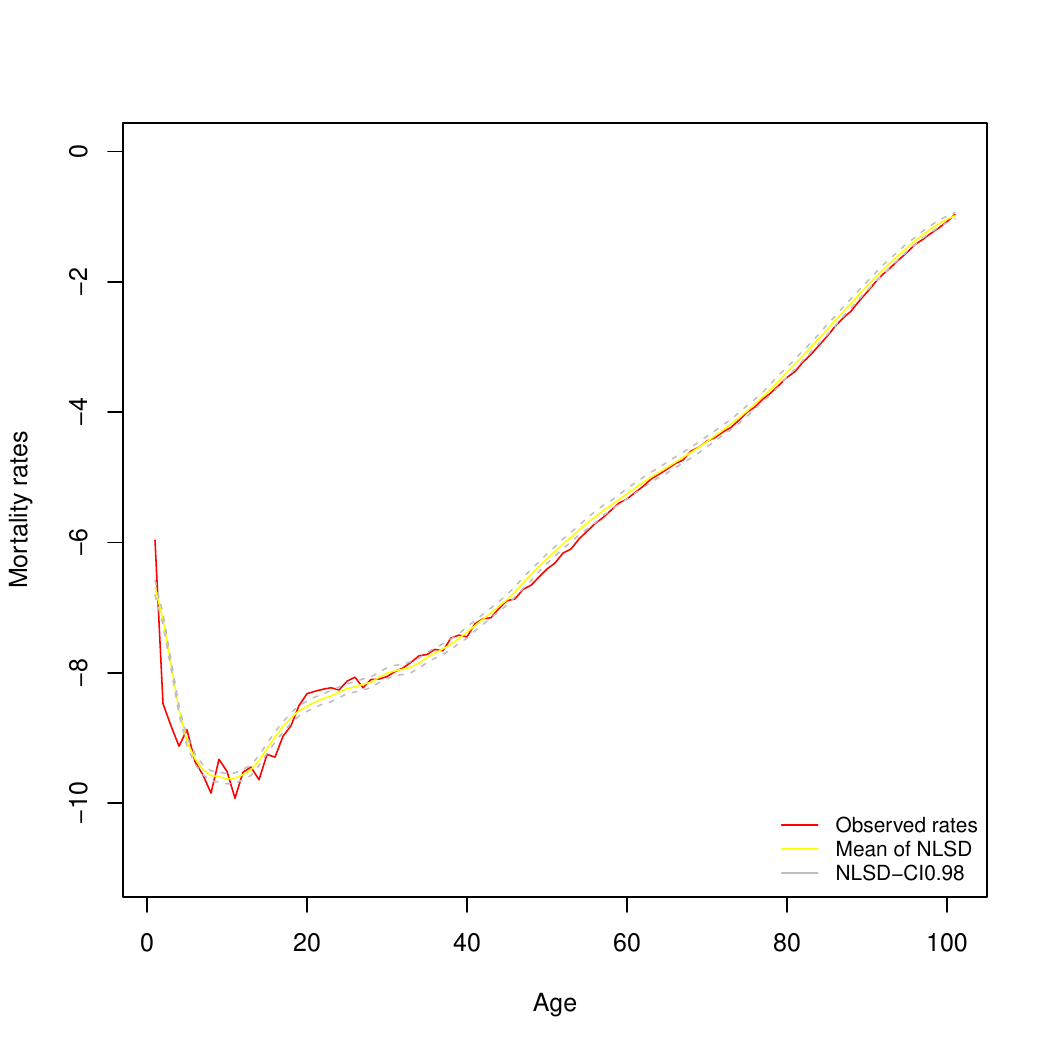}
        \caption{{Confidence Level: $1-\alpha=0.98$}}
        \label{Fig24}
    \end{subfigure}

    \vspace{0.5cm} 

    \begin{subfigure}[b]{0.45\textwidth}
        \centering
        \includegraphics[scale=0.5,width=\textwidth]{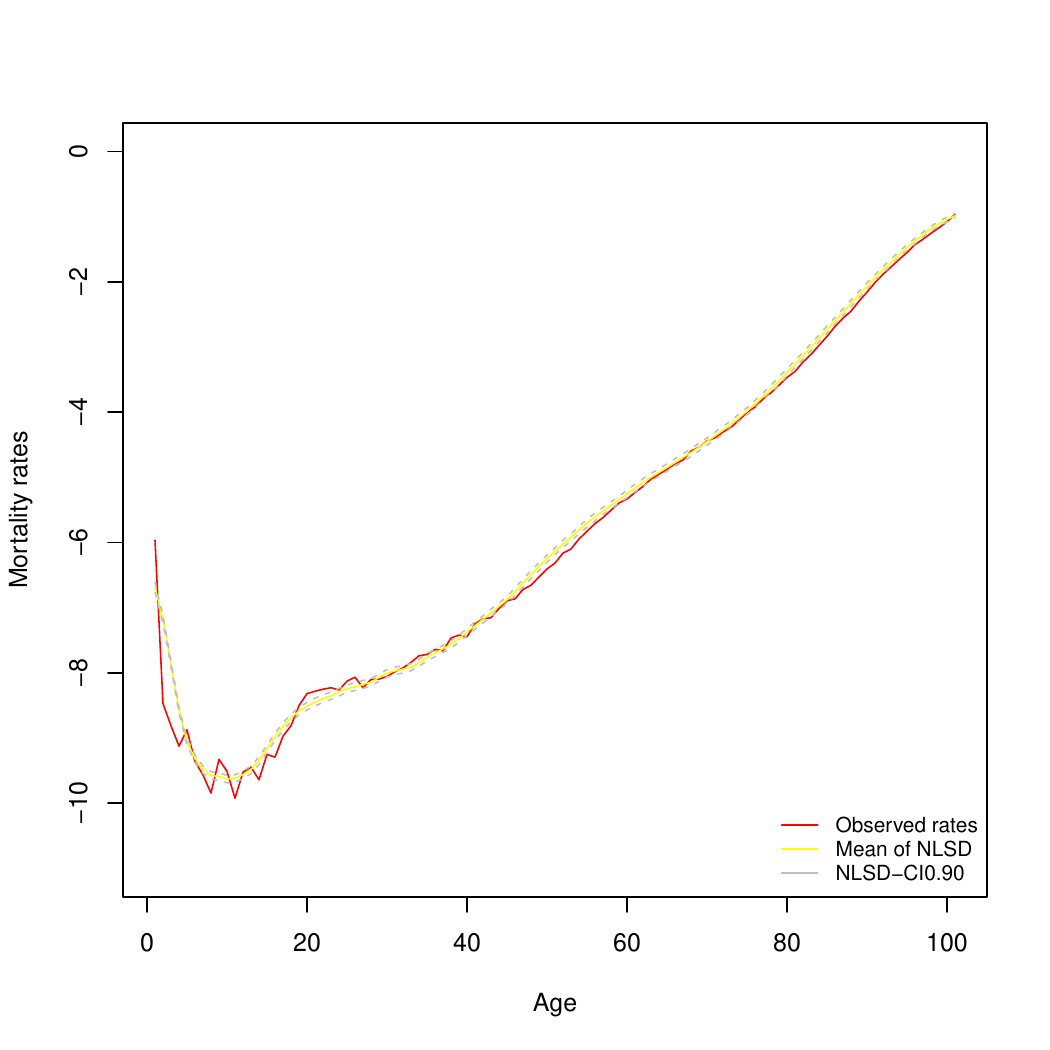}
        \caption{{Confidence Level: $1-\alpha=0.9$}}
        \label{Fig25}
    \end{subfigure}
    \hfill
    \begin{subfigure}[b]{0.45\textwidth}
        \centering
        \includegraphics[scale=0.5,width=\textwidth]{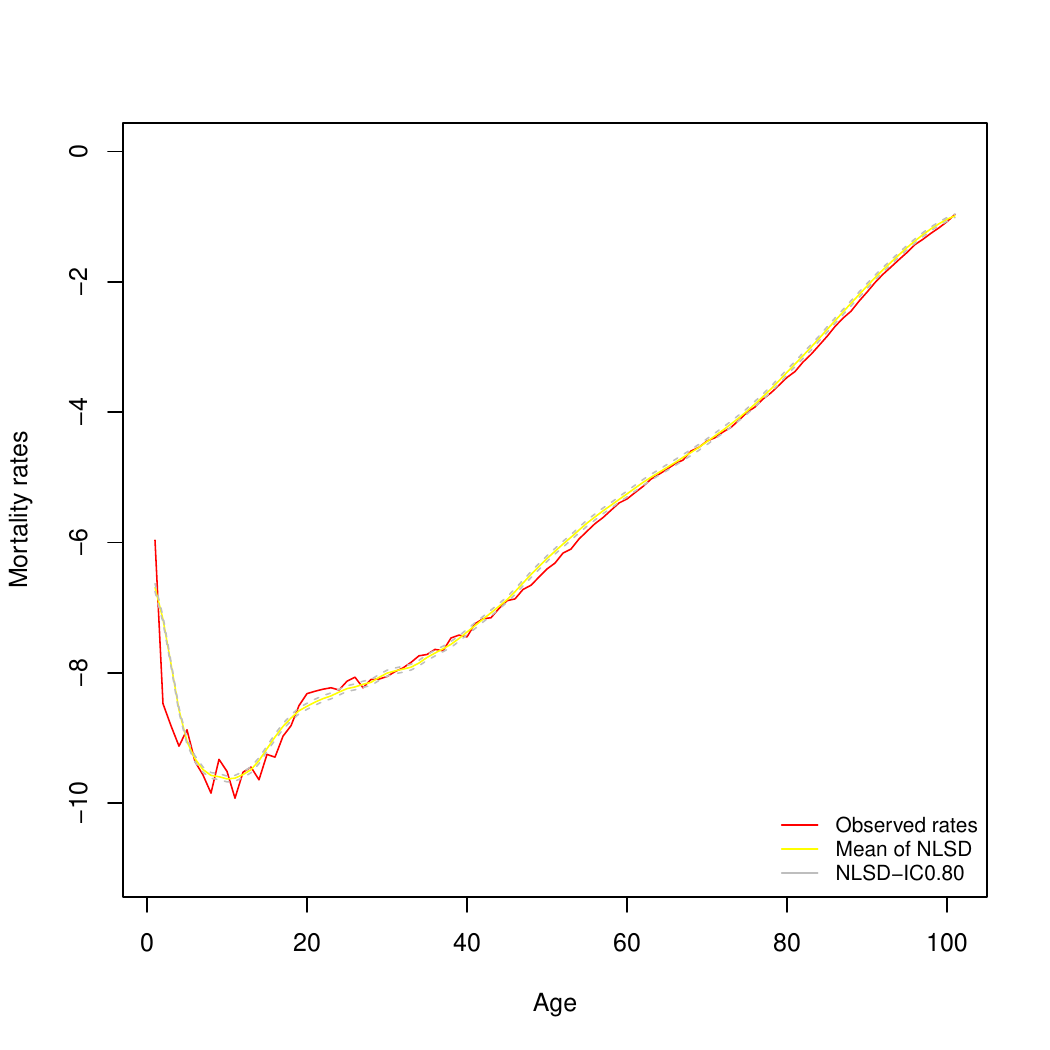}
        \caption{{Confidence Level:$1-\alpha=0.8$}}
        \label{Fig26}
    \end{subfigure}

    \caption{Confidence Intervals for several confidence levels ($b=0.025$).}
    \label{fig:figura_CIb025}
\end{figure}

Figures \ref{Fig21}-\ref{Fig23} and Figures \ref{Fig24}-\ref{Fig26} show, for the 5-year horizon of forecasting,
if the observed rates belong or not to the intervals $IC_{0.98}$, $IC_{0.90}$
and $IC_{0.80}$ for $b=0.1$ and $b=0.025$.

\begin{figure}[htbp]
\centerline{{\includegraphics[scale=0.5,angle=0]{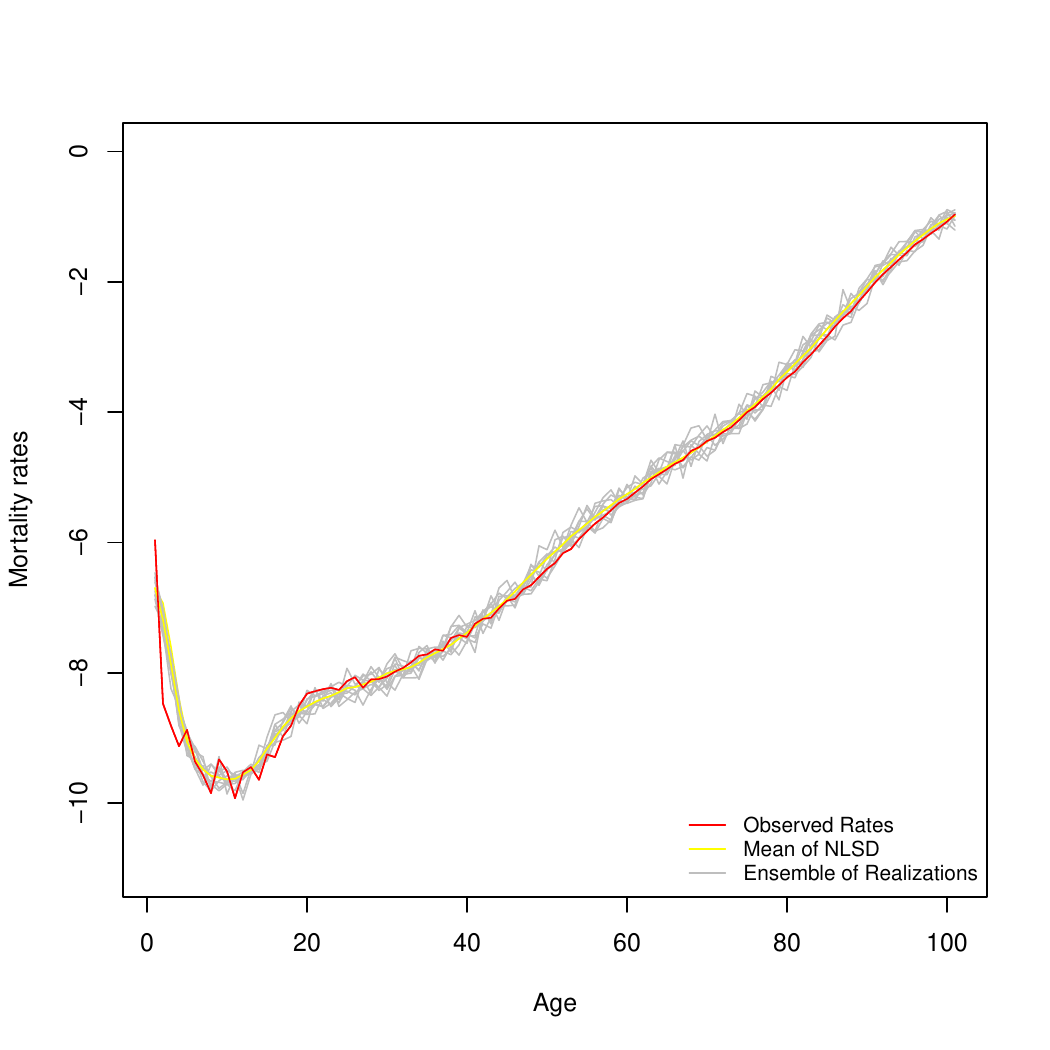}}}\caption{{\protect\footnotesize {Ensemble
of realizations: $b=0.1$}}}%
\label{Fig3}%
\end{figure}

In the context of real applications, forecasting several plausible scenarios
often requires more than just the mean trajectory. For example, when the
mortality rates are used as input data in nonlinear estimations, as in the
calculation of the cost of claims using mechanisms based on compound interest
rates, it becomes convenient to account for random fluctuations. As we can see
in Figure \ref{Fig3}, the proposed method allows us to estimate not only the
mean trajectories but also an arbitrary number of equally probable realizations.

Further, with the aim of comparison, we apply different methods using the same
data in the period 1908-2018 and forecast the mortality rates for the period
2019-2023. Then, we calculate the indicators which have been defined previosly.

As we said before, initially we have considered several methods such as the
Lee-Carter (LC), Renshaw-Haberman (RH) and CBD methods, and the models M3-M7.
With the aim of facilitating the interpretation of the results, we have
selected the best-fitting methods. The selection is due using the Akaike
Information Criterion (AIC) and the Bayesian Information Criterion (BIC). With
these indicators, and using the package StMoMo in the R-software, it has been
determined that the Renshaw-Haberman metod and the Lee-Carter method are the
most suitable.

\begin{figure}[htbp]
\centerline{{\includegraphics[scale=0.5,angle=0]{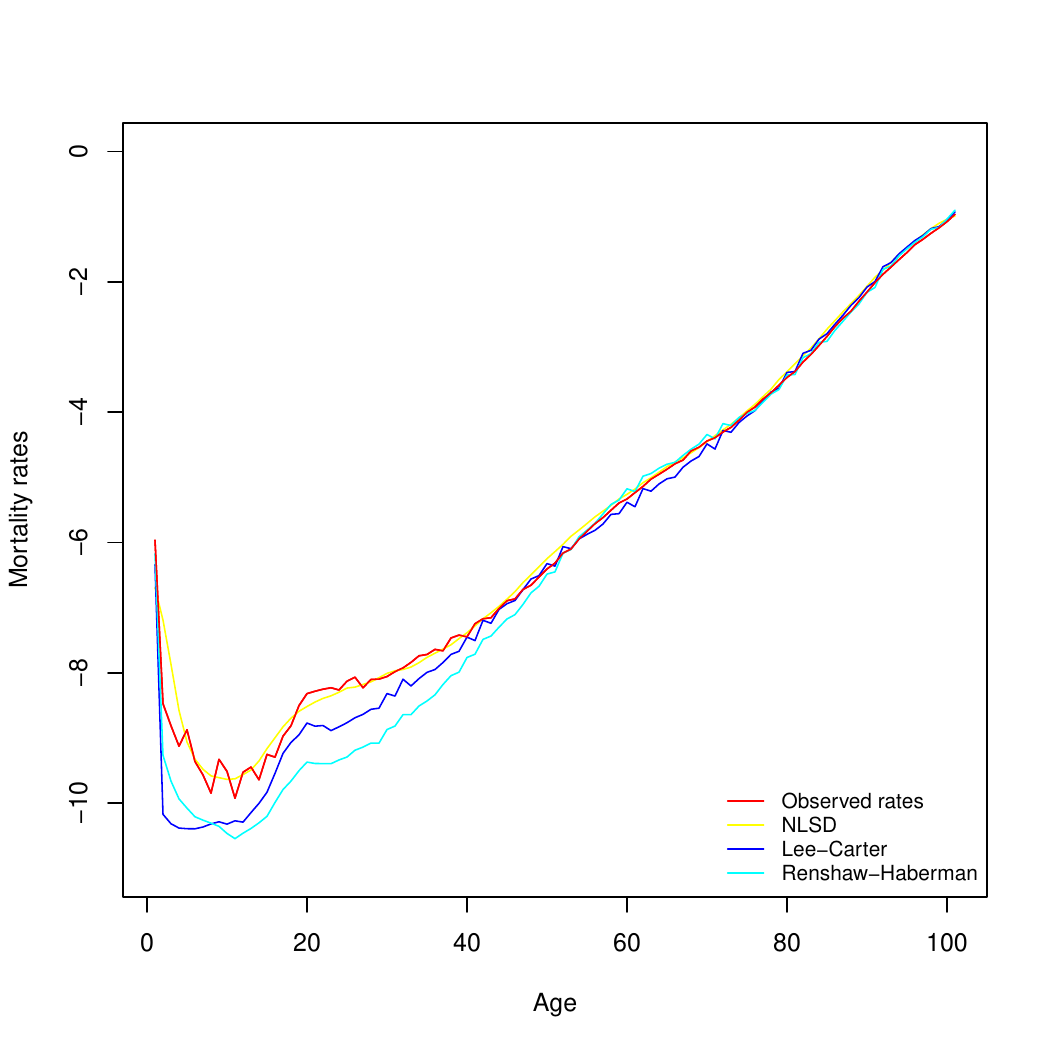}}}\caption{{\protect\footnotesize {Mean
trajectories}}}%
\label{Fig4}%
\end{figure}

Figure \ref{Fig4} shows us the results for the year 2023: the observed data
and the mean value of the trajectories for each technique. In the NLSD method,
$b=0.1$.

Figure \ref{Fig4} provides evidence about the suitability of the proposed
method. Even though the proposed technique seems to be, graphically, the
better technique, it is important to note that none of the three evaluated
techniques has been specifically calibrated, and the default parameter values
of the StMoMo package have been used for the LC and RH methods.

Complementarialy, we can use the quantitave measures to determine the goodness
of each model and to compare them.

We start with the count indicators. Table 1 shows the number of ages (for each
year into the period of the validation) that do not belong to the confidence
interval, $IC_{0.98}$, and for each of the evaluated methods. Tables 2 and 3
show the same information for the confidence intervals $IC_{0.90}$ and
$IC_{0.80}$. Our method has been tested with different values of the parameter
$b$ ($0$.$1,\ 0$.$05$ and $0$.$025$), which determines the intensity of the noise.

\begin{center}
$\underset{\text{{\large Table 1: }}I_{c,0.98}^{t}}{%
\begin{tabular}
[c]{|r|c|c|c|c|c|}\hline
$\text{Method%
$\backslash$%
year}$ & 2019 & 2020 & 2021 & 2022 & 2023\\\hline
NLSD\ 0.1 & \multicolumn{1}{|r|}{7} & \multicolumn{1}{|r|}{6} &
\multicolumn{1}{|r|}{4} & \multicolumn{1}{|r|}{8} & \multicolumn{1}{|r|}{5}%
\\\hline
NLSD\ 0.05 & \multicolumn{1}{|r|}{17} & \multicolumn{1}{|r|}{56} &
\multicolumn{1}{|r|}{21} & \multicolumn{1}{|r|}{18} & \multicolumn{1}{|r|}{33}%
\\\hline
NLSD\ 0.025 & \multicolumn{1}{|r|}{50} & \multicolumn{1}{|r|}{75} &
\multicolumn{1}{|r|}{53} & \multicolumn{1}{|r|}{39} & \multicolumn{1}{|r|}{67}%
\\\hline
LC & \multicolumn{1}{|r|}{65} & \multicolumn{1}{|r|}{80} &
\multicolumn{1}{|r|}{70} & \multicolumn{1}{|r|}{63} & \multicolumn{1}{|r|}{71}%
\\\hline
RH & \multicolumn{1}{|r|}{71} & \multicolumn{1}{|r|}{88} &
\multicolumn{1}{|r|}{81} & \multicolumn{1}{|r|}{81} & \multicolumn{1}{|r|}{84}%
\\\hline
\end{tabular}
}$

\medskip

$\underset{\text{{\large Table 2:\ }}I_{c,0.9}^{t}}{%
\begin{tabular}
[c]{|r|r|r|r|r|r|}\hline
$\text{Method%
$\backslash$%
year}$ & 2019 & 2020 & 2021 & 2022 & 2023\\\hline
NLSD\ 0.1 & 10 & 18 & 9 & 11 & 16\\\hline
NLSD\ 0.05 & 26 & 66 & 36 & 28 & 59\\\hline
NLSD\ 0.025 & 64 & 84 & 66 & 50 & 82\\\hline
LC & 77 & 93 & 80 & 72 & 82\\\hline
RH & 78 & 93 & 82 & 87 & 87\\\hline
\end{tabular}
}$

\medskip

$\underset{\text{{\large Table 3:\ }}I_{c,0.8}^{t}}{%
\begin{tabular}
[c]{|r|r|r|r|r|r|}\hline
$\text{Method%
$\backslash$%
year}$ & 2019 & 2020 & 2021 & 2022 & 2023\\\hline
NLSD\ 0.1 & 15 & 44 & 17 & 15 & 28\\\hline
NLSD\ 0.05 & 41 & 72 & 48 & 34 & 64\\\hline
NLSD\ 0.025 & 72 & 91 & 75 & 67 & 86\\\hline
LC & 79 & 95 & 85 & 75 & 88\\\hline
RH & 85 & 94 & 89 & 91 & 91\\\hline
\end{tabular}
}$

\medskip
\end{center}

By analyzing the tables of the count-based indicators, we can conclude that
the proposed technique captures the observed mortality more effectively. This
suggests that, in this regard, it provides a better fit than the other
approaches. However, this does not necessarily imply that the technique is
more accurate, as the result may be due to a more conservative forecast, so we
have to consider other indicators as well.

Complementary to the count-based indicators, error metrics may be useful for
assessing whether the technique is as accurate as, or more or less accurate
than, the LC and RH models. Using the indicators (\ref{error_measures}) over
the full 0-100 age range we obtain the results in Tables 4 and 5.

\begin{center}
$\underset{\text{{\large Table 4:\ }}I_{MqD}^{t}}{%
\begin{tabular}
[c]{|r|l|l|l|l|l|}\hline
$\text{Method%
$\backslash$%
year}$ & 2019 & 2020 & 2021 & 2022 & 2023\\\hline
NLSD\ 0.1 & 2.256290e-05 & 2.236960e-05 & 2.260896e-05 & 2.961290e-05 &
\textbf{1.518093e-05}\\\hline
NLSD\ 0.05 & 1.103389e-04 & 1.091354e-04 & 1.090309e-04 & 1.070186e-04 &
1.554116e-04\\\hline
NLSD\ 0.025 & \textbf{1.075917e-05} & \textbf{1.107349e-05} &
\textbf{1.083946e-05} & \textbf{1.529299e-05} & 1.772229e-05\\\hline
LC & 3.135750e-05 & 3.088876e-05 & 3.099892e-05 & 2.791308e-05 &
4.069674e-05\\\hline
RH & 3.186545e-05 & 3.139071e-05 & 3.097631e-05 & 2.684819e-05 &
1.960841e-05\\\hline
\end{tabular}
}$

\medskip

$\underset{\text{{\large Table 5:\ }}I_{MRqD}^{t}}{%
\begin{tabular}
[c]{|r|l|l|l|l|l|}\hline
$\text{Method%
$\backslash$%
year}$ & 2019 & 2020 & 2021 & 2022 & 2023\\\hline
NLSD\ 0.1 & 1.162379e-04 & 1.168913e-04 & 1.174218e-04 & 1.815611e-04 &
9.268740e-05\\\hline
NLSD\ 0.05 & 4.740293e-04 & 4.682191e-04 & 4.692663e-04 & 4.908082e-04 &
7.605474e-04\\\hline
NLSD\ 0.025 & \textbf{7.714122e-05} & \textbf{7.831104e-05} &
\textbf{7.724071e-05} & \textbf{1.283547e-04} & 1.373643e-04\\\hline
LC & 1.075780e-04 & 1.061766e-04 & 1.066380e-04 & 1.506715e-04 &
2.467346e-04\\\hline
RH & 2.252991e-04 & 2.267182e-04 & 2.237182e-04 & 1.858389e-04 &
\textbf{1.321769e-04}\\\hline
\end{tabular}
}$
\end{center}

Tables 4 and 5 allow us to compare the methods. We have highlighted in bold
the values with lowest error for each year. According to Table 4, the method
NLSD has the lowest error accross all the years. In particular, NLSD with
$b=0.025$ outperforms both the LC and RH methods in all the years. In Table 5,
NLSD (with $b=0.025$) is the most accurate method in four out of five years,
while RH performs best in one (although the error value for NLSD is nearly
identical to that of RH in that case). The LC model does not achieve the
lowest error in any of the five years for either indicator.

Central measures are also very useful to assess the accuracy of the methods.
In Table 6 one can see the values of the indicators (\ref{Central1}) in the
year 2023.

\begin{center}

$\underset{\text{{\large Table 6: }}I_{CT,\tau}^{t}}{%
\begin{tabular}
[c]{|r|r|r|}\hline
& $I_{CT,1}^{2023}$ & $I_{CT,2}^{2023}$\\\hline
NLSD\ 0.1 & 0.18 & 149.54\\\hline
NLSD\ 0.05 & 0.36 & 623.98\\\hline
NLSD\ 0.025 & 0.70 & \textit{2490.81}\\\hline
LC & 0.88 & 2109.89\\\hline
RH & 1.74 & 103496.19\\\hline
\end{tabular}
}$

\end{center}

The indicator $I_{CT,1}^{2023}$ yields better results for NLSD across all
three levels of noise intensity. The indicator $I_{CT,2}^{2023}$ is better for
NLSD for $b=1$ and $b=0.05$, while LC slightly outperforms NLSD for $b=0.025$.
The values obtained for RH are significantly worse than for the other methods.

\begin{figure}[htbp]
\centerline{{\includegraphics[scale=0.4,angle=0]{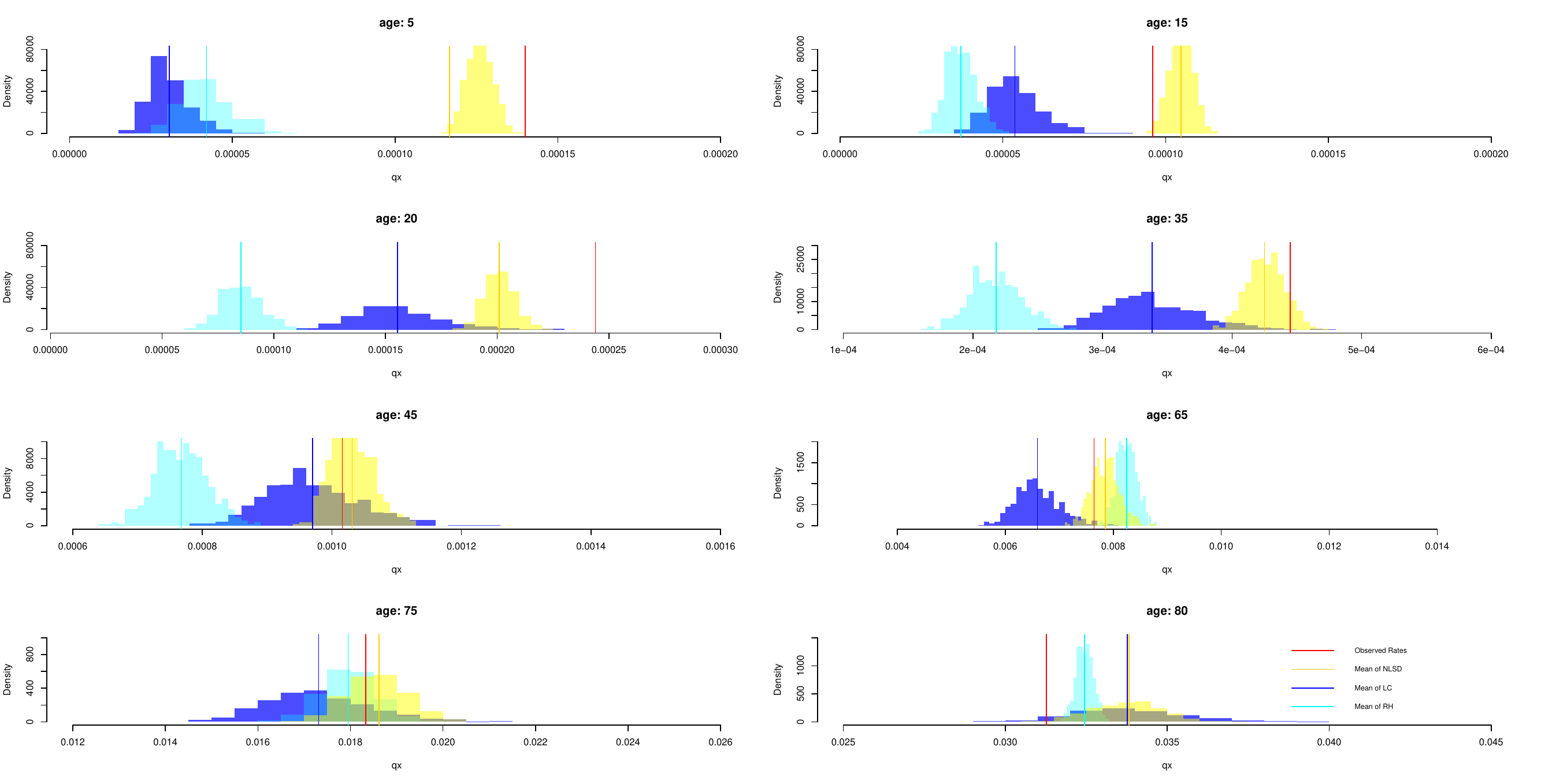}}}\caption{{\protect\footnotesize {Density function, $b=0.025$}}}%
\label{FigDensity}%
\end{figure}

Figure \ref{FigDensity} shows the estimated density functions for several ages and each
method. From this picture we can see graphically the variability, which is
different for each age and method. This change in variability from one
technique to another, when directly comparing the estimated mean values and
the observed values, highlights the importance of accounting for such
variability in order to accurately assess the precision of each technique.

When $b=0.025$, for most ages, the NLSD method yields mean values closest to the observed values, while maintaining a level of variability comparable to the other methods. At the oldest age (80), the LC method provides the best coverage, albeit at the cost of high variability in the realizations. Moreover, for younger and middle ages, the realizations from the LC method fail to cover the observed value, despite exhibiting greater variability.

\bigskip

\begin{figure}[htbp]
\centerline{{\includegraphics[scale=0.5,angle=0]{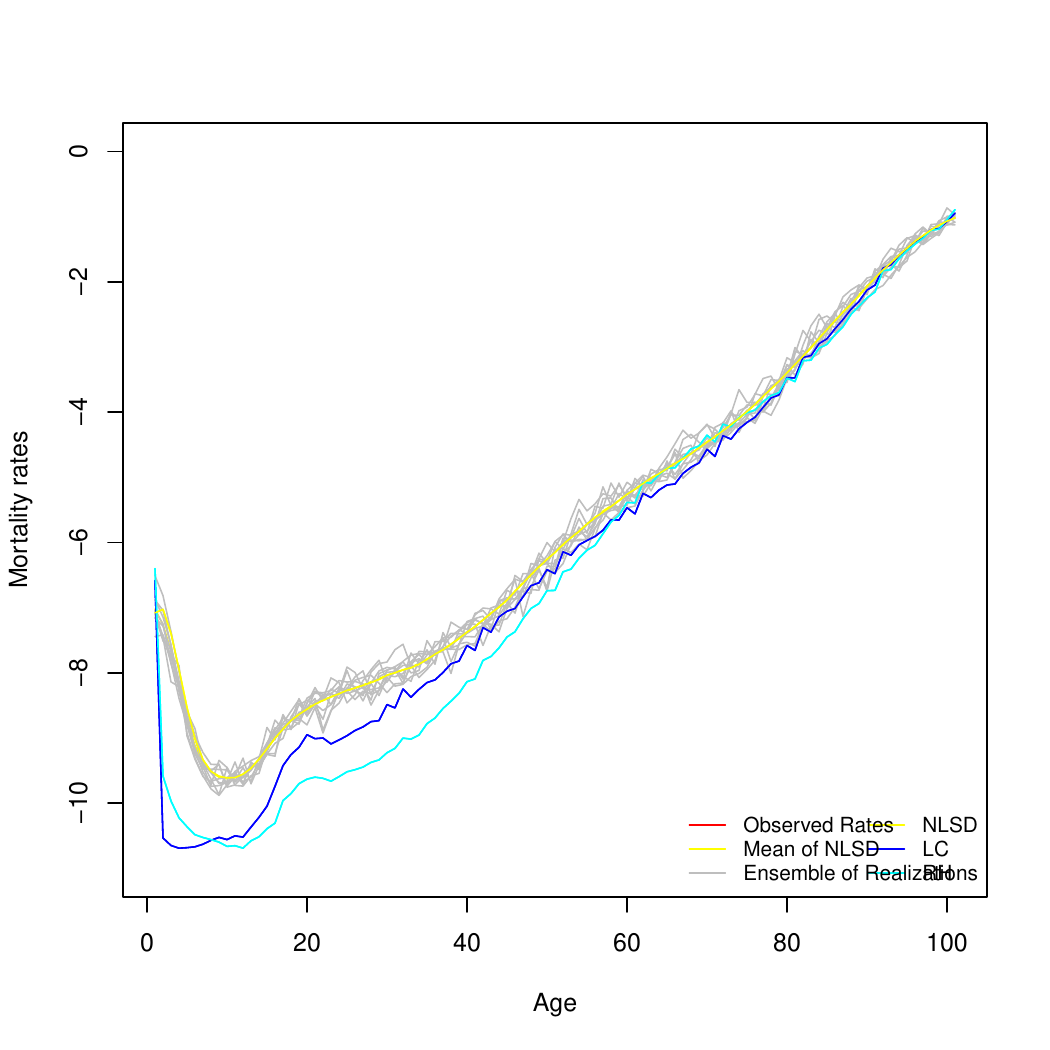}}}\caption{{\protect\footnotesize {Mean
trajectories}}}%
\label{Fig5}%
\end{figure}

We have seen that the NLSD method can be applied to real-world scenarios with
high short-term accuracy. To assess whether this technique is also effective
in the medium or long term, we examine whether its predicted values align with
those obtained from the RH and LC methods. Figure \ref{Fig5} displays the mean
prediction trajectories over a 10-year horizon (with predictions for 2028
based on observed data up to 2018).

We observe that the predicted values are similar in magnitude across the three
techniques. However, there exist differences at the youngest and oldest ages.
The predictions diverge most at the earliest ages, where the NLSD model
produces the highest values, followed by LC with intermediate values, and RH
with the lowest.

From a qualitative point of view, it is worth noting that the analysis of
historical time series indicates a decreasing intensity of the social hump
over time. In this regard, the NLSD model exhibits a more realistic behavior
compared to the LC and RH models.

It is also important to highlight that the LC and RH models generate
predictions using autoregressive and moving average time series models
(ARIMA), which possess the following characteristics:

\begin{enumerate}
\item They make linear forecasts by extrapolating the dynamics of the most
recent values. For instance, recent improvement in mortality rates may be
projected forward, potentially leading to underestimations of future mortality levels.

\item They exhibit an \textit{uncertainty cone} that grows rapidly so their
predictive performance deteriorates significantly as the forecast horizon
increases. In Economics, for example, it is common to use a 12-step monthly
forecast horizon. For annual forecasts, typical horizons are 3, 5, or 7 years.
\end{enumerate}

\begin{figure}[htbp]
\centerline{{\includegraphics[scale=0.5,angle=0]{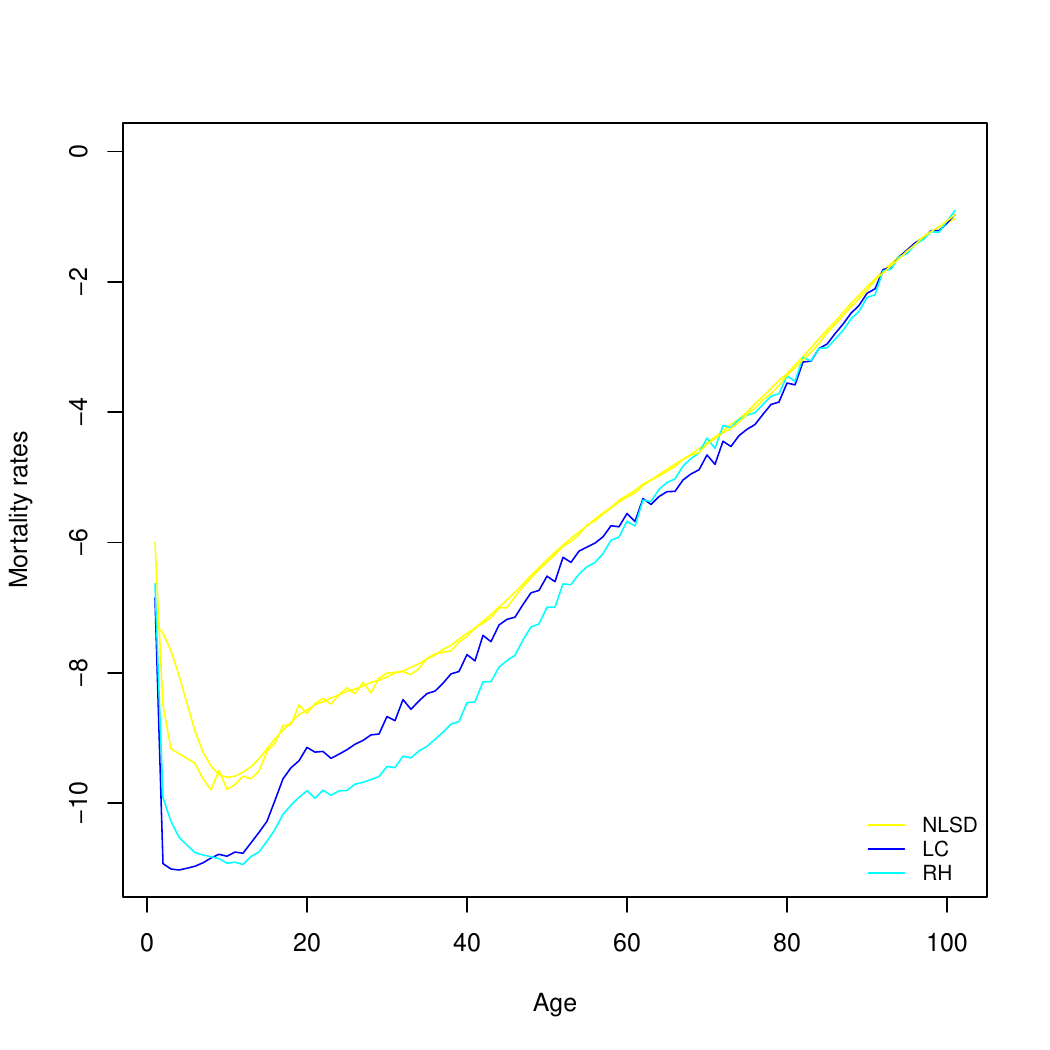}}}\caption{{\protect\footnotesize {Mean
trajectories}}}%
\label{Fig6}%
\end{figure}

Figure \ref{Fig6} shows the predictions for the year 2033 (15-year horizon).
This figure exhibits a similar pattern to that observed in the 10-year
horizon, but the differences between the methods become more pronounced.

\section{Conclusions}

This work proposes a method for modeling and forecasting mortality rates. It
constitutes an improvement over previous studies (\cite{MoVa14},
\cite{MoVa01}, \cite{CarMoVa02}), by incorporating both the historical
evolution of the mortality phenomenon and its random behavior.

The first part of the study introduces the NLSD model and analyzes
mathematical properties such as the existence of solutions and their
asymptotic behavior. The second part presents an application of the NLSD
model. For this purpose, the Euler--Maruyama method is applied to data
obtained from the Human Mortality Database \cite{HMD}. The choice of the HMD
is justified by the fact that it contains mortality datasets from a large
number of countries, all of which have been methodologically harmonized. The
use of Spanish data is arbitrary; the method has also been tested with data
from other countries, such as the UK, although we do not show these results in
this paper.

To assess the validity of the proposed method, the observation period was
divided into two subsets: one for fitting and calibration ($1908-2018$), and
another for validation ($2019-2023$).

The evaluation was carried out by comparing the proposed model with other
widely used approaches, such as the LC, RH, CBD, and M3--M7 models.
Count-based, error-based and central metrics were used in the comparison. The
NLSD model achieved the best results for all years within the validation period.

Based on this study, we can conclude that the NLSD model should be regarded as
a promising alternative to classical models. While a more exhaustive
validation remains to be conducted, the method has shown the best performance
among the models tested.

As extensions of this study, we propose conducting a sensitivity analysis of
the parameters, as well as an exhaustive comparison across different regions
and time periods. From a technical perspective, it would be valuable to
incorporate optimization techniques for parameter estimation and to assess the
applicability of cross-validation strategies.

\bigskip

\textbf{Acknowledgements.} The research has been partially supported by the
Spanish Ministerio de Ciencia e Innovaci\'{o}n (MCI), Agencia Estatal de
Investigaci\'{o}n (AEI) and Fondo Europeo de Desarrollo Regional (FEDER) under
the project PID2021-122991NB-C21.

\textbf{Conflict of interest.}
The authors declare that they do not have any conflict of interest.


\end{document}